\documentclass{amsart}
\usepackage{amsmath, amssymb, amsthm, epsf, graphicx, amscd, amsfonts, amscd, multirow}
\usepackage[all]{xy}
\usepackage{setspace}
\usepackage{xcolor}

\newtheorem{theorem}{Theorem}[section]
\newtheorem{lemma}[theorem]{Lemma}

\newtheorem{remark}[theorem]{Remark}
\newtheorem{corollary}[theorem]{Corollary}

\def\bb #1{ {\mathbb #1} }
\def\c #1{ {\mathcal #1} }

\def\b #1{ {\bf #1} }

\numberwithin{equation}{section}

\begin{document}
\title{Maximal Unipotent Monodromy and $p$-Adic Eigenvalue Formulas}

\author{Alan Adolphson}
\address{Department of Mathematics, Oklahoma State University, Stillwater, Oklahoma 74078}
\email{alan.adolphson@okstate.edu}
\author{Steven Sperber} 
\address{School of Mathematics, University of Minnesota, Minneapolis, Minnesota 55455}
\email{sperber@math.umn.edu}
\author{Hao Zhang}
\address{School of Mathematics and Shing-Tung Yau Center, Southeast University, Nanjing 211189, People’s
Republic of China}
\email{hzhang21@seu.edu.cn}

\date{\today}
\begin{abstract}
We establish $p$-adic formulas for the eigenvalues of Frobenius in the case of the family of hyper\-kloosterman sums. The deformation equation in this case has maximal unipotent monodromy, and the interplay of Frobenius and the differential equation plays an important role in establishing $p$-adic analytic formulas for the  eigenvalues. 
\end{abstract}
\maketitle

\section{Introduction}

We describe here some of the key properties we use. We consider the hyperkloosterman family of exponential sums. Let $\bb F_q$ be a finite field of characteristic $p$, with $q = p^a$. Let
\[
K_n( \Lambda, x) := x_1 + \cdots + x_n + \frac{ \Lambda}{x_1 \cdots x_n} \in \bb Z[\Lambda, x_1^{\pm 1},\dots,x_n^{\pm 1}].
\]
Let $\Theta$ be a non-trivial additive character of  ${\bb F}_p$
and for any positive integer $j$ let $\Theta_j = \Theta \circ \text{Tr}_{\mathbb{F}_{q^j}/\mathbb{F}_p}$ be the additive character of ${\bb F}_{q^j}$ defined by taking the composition of $\Theta$ with the trace map. 
For $\bar{\lambda}\in{\bb F}_q$, define
\[ Kl_n(\bar{\lambda})/{\bb F}_{q^j} = \sum_{x=(x_1,\dots,x_n)\in({\bb F}_{q^j}^\times)^n} \Theta_j(K_n(\bar{\lambda},x)). \] 

It is well-known that for $\bar{\lambda}\in{\bb F}_q^\times$ the associated $L$-function
\[
L(Kl_{n}/\mathbb{F}_q, \bar \lambda,  T)^{(-1)^{n+1}} = \exp \ \bigg(\sum_{j=1}^{\infty} \frac{(Kl_{n}(\bar \lambda)/\mathbb{F}_{q^j}) \ T^j}{j}\bigg)
\]
is a polynomial in $T$ of degree $n+1$ (see \cite{AS0}), with reciprocal roots $ \{\alpha_i \}_{i=0}^n$ with 
\begin{equation}
\alpha_i = q^i u_i(\bar{\lambda})                 
\end{equation}
where   $u_i(\bar{\lambda})$ is $p$-adically a principal unit. Our goal in the present work is to give a $p$-adic formula which realizes $u_i(\bar{\lambda})$  for each $i=0,1,\dots,n$ as a special value of a suitable $p$-adic analytic function at $\hat{\lambda}$, the Teichm\"uller lifting of $\bar{\lambda}.$ Such formulas have appeared in the literature for the unique unit root $u_0(\bar{\lambda})$ in the hyperkloosterman family (and other families when such a unique unit root exists). In the generalized Airy case (Zhang \cite{Z}), formulas were proven for all reciprocal roots and the methods used there apply here as well. The approach in \cite{Z} looks at the Frobenius action on the exterior powers of the deformation equation. A related approach developed by Dwork \cite{D1,D2} involves a factorization of normalized local solution matrices to the deformation equation using a global factor. 
 
In \cite{AS} Adolphson and Sperber show that equation (1.1) holds in all characteristics. That work required us to work with the splitting function $\theta_{\infty}$.  We will prefer here to use the $\theta_1$ splitting function of Dwork.  (The loss that this choice incurs in the restriction of the characteristic $p$ is necessitated by our need for a compatible  $p$-adic formulation of the functional equation.)  This necessitates that we restrict to the case $p>n+1$ and this hypothesis will hold throughout the paper. 

As usual, let $\mathbb{Q}_p$ be the field of $p$-adic numbers, $\mathbb{Z}_p$ be the ring of $p$-adic integers,  $\mathbb{Q}_q$ be the unramified extension of $\mathbb{Q}_p$ having residue class field $\mathbb{F}_q$, the finite field of $q$ $(= p^a)$ elements.  Let $\Omega_1 = \mathbb{Q}_p(\zeta_p)$, and $\Omega = \mathbb{Q}_q(\zeta_p)$ be the respective cyclotomic extensions of $\mathbb{Q}_p$ and $\mathbb{Q}_q$ and $\c O$ the ring of integers of  $\Omega$. The nontrivial additive character $\Theta$ of ${\mathbb F}_p$ determines a uniformizer $\pi$ for $\Omega$ by the congruence $\Theta(1)\equiv 1+\pi\pmod{\pi^2}$.  We define a ring of functions over which we will work. As in previous works (Haessig-Sperber (\cite{HS}), we define a weight function $w$ with $w(\lambda) = n+1$ and for any $\alpha \in \mathbb{Z}_{\geq 0}$ we have $w(\lambda^{\alpha}) (= w(\alpha)) = \alpha(n+1)$. We will work over the ring of functions (where $\lambda$ is a variable) 
\[ \c{C}_0 = \bigg\{ \sum_{\alpha \in \mathbb{Z}_{\geq 0}} c(\alpha) \pi^{w(\alpha)} \lambda^{\alpha}\mid c(\alpha) \in \c{O}, c(\alpha) \rightarrow 0 \  \text{as} \  \alpha \rightarrow \infty \bigg\}.  \]
In addition, for future use we also define a  $\c{C}_0$-module of ``bounded functions"  
\[ \tilde{\c C_0} =  \bigg\{ \sum_{\alpha \in \mathbb{Z}_{\geq 0}} c(\alpha) \pi^{w(\alpha)} \lambda^{\alpha}\mid c(\alpha) \in \c O \bigg\}.  \]
Of course, $ \c{C}_0 \subset  \tilde{\c C_0}$  and we give both of these a norm by defining for $\xi =  \sum_{\alpha \in \mathbb{Z}_{\geq 0}}  c(\alpha) \pi^{w(\alpha)} \lambda^{\alpha} \in  \tilde{\c C_0}, \  |\xi| = \ \sup_{\alpha \in \mathbb{Z}_{\geq 0}} \  |c(\alpha)| $. 

Let $\c M$ denote the monoid generated by $\langle x_1,\dots,x_n,\frac{\lambda}{x_1x_2\cdots x_n}\rangle$.  
We recall a weight function can be defined on elements of $ \lambda^{\sigma}x^u \in \c M$ via $W(\sigma, u) = \sigma(n+1) + \sum_{i=1}^n u_i$. 
If we set $m(u) = \max\{0,-u_1,\dots,-u_n\}$ for any  $u\in{\mathbb Z}^n$, then the usual polyhedral weight of $u$
is $w(u) = W(m(u), u) = (n+1)m(u) + \sum_{i=1}^n u_i$. We may define closed sectors $\{I_j\}_{0 \leq j \leq n}$ of ${\mathbb{Z}^n}$ as follows: $I_0 = (\mathbb{Z}_{\geq 0})^n $ and for $1 \leq j \leq n$, $I_j =  \{ u \in \mathbb{Z}^n \mid -u_j = m(u) \}$.
With these definitions, $W(\sigma, u)$ is linear with respect to $\sigma$ and $u$ while $w(u+v)\le w(u)+w(v)$ with equality holding if and only if $u$ and $v$ belong to the same closed sectors.

We define a ring of functions 
$\c L$, which is a $\c C_0$ algebra, via
\[  \c L = \bigg\{\sum _{\mu \in \mathbb{Z}^n} C(\mu) \lambda^{m(\mu)} \pi^{w(\mu)} x^{\mu} \mid C(\mu) \in \c C_0, C(\mu) \rightarrow 0 \ \text{as} \ \mu \rightarrow \infty \bigg \}.
\]
Proceeding then $p$-adically, we define a complex $\Omega^{\bullet}_{\c L} $ of $\c L$-modules as follows.  For each $\b{j}$ an ordered subset of $\c S = \{1, 2,\dots,n\}$ of cardinality $i$ such that $\b j= \{1 \leq j_1 < j_2 < \cdots <j_i \leq n\}$, we set $\Omega^i_{\b j} =  \c L \frac{dx_{\b j}}{{x_{\b j}}}$ where $ \frac{dx_{\b j}}{{x_{\b j}}}=\frac{dx_{j_1}}{x_{j_1}} \wedge\cdots\wedge \frac{dx_{j_i}}{x_{j_i}}$, and set $\Omega^i = \bigoplus \Omega^i_{\b j}$ the sum taken over all ordered subsets $\b j$ of $\c S$ of cardinality $i$. We define the boundary operator as follows.  For $\xi = \alpha \frac{dx_{\b j}}{x_{\b j}}$, we set $\nabla \xi  = \sum_{\ell =1}^n D_{\ell} (\alpha) \frac{dx_{\ell}}{x_{\ell}}  \wedge \frac{dx_{\b j}}{x_{\b j}}$. Here $D_{\ell}$ is the twisted differential operator 
\[ 
D_{\ell} = \exp( -\pi K_n(\Lambda,x)) \circ x_{\ell}\frac{{\partial}}{\partial x_{\ell}} \circ \exp (\pi K_n(\Lambda, x)) = E_{\ell} + H_{\ell},
\]
where we use the notation $E_{\ell} = x_\ell\frac{\partial}{\partial x_{\ell}}$ and $H_{\ell} = E_{\ell} (\pi K(\Lambda. x)) =  \pi (x_{\ell} - \frac{\Lambda}{\prod _{i=1}^n x_i} )$.
It is a well-known result \cite{S} that the cohomology of this complex vanishes except in top dimension and that $H^n(\Omega^{\bullet})$ is a free $\c{C}_0$-module of rank $n+1$ with a basis given by $ \{ \epsilon_i=\pi^ix_1\cdots x_i \}_{i=0}^n$. We use the convention $\epsilon_0 =1$.

Following Dwork, we define a dual space,
$\c L^{*}  = \bigg\{\sum _{\mu \in \mathbb{Z}^n} C^{*}(\mu) \lambda^{-m(\mu)} \pi^{-w(\mu)} x^{-\mu} \mid C^{*}(\mu) \in \c C_0 \bigg\}$, a $\c C_0$ -module.  There is  a pairing $\langle \c L, \c L^{*} \rangle \rightarrow \c C_0$, where  for $\eta =  \sum _{\mu \in \mathbb{Z}^n} \eta(\mu) \lambda^{m(\mu)} \pi^{w(\mu)} x^{\mu} \in \c L$ with  $\eta(\mu) \in \c C_0$, $\eta(\mu) \rightarrow 0 \ \text{as} \ \mu \rightarrow \infty$ and $\eta^{*} =  \sum _{\mu \in \mathbb{Z}^n} \eta^{*}(\mu) \lambda^{m(\mu)} \pi^{w(\mu)} x^{\mu}\in\c L^*$ with $\eta^{*}(\mu) \in \c C_0$, then $\langle \eta, \eta^{*} \rangle = \sum _{\mu \in \mathbb{Z}^n} \eta(\mu) \eta^{*}(-\mu)  \in \c C_0$. In an altogether analogous way, we construct a dual complex $(\Omega_{\c L^{*}})^{\bullet}$  working with $\c L^{*}$ as coefficients and dual differential operators 
\[ 
D_{\ell}^{*} = \exp( \pi K_n(\Lambda,x)) \circ -x_{\ell}\frac{{\partial}}{\partial x_{\ell}} \circ \exp (-\pi K_n(\Lambda, x)) = -E_{\ell} + H_{\ell}.
\]

Fix the basis $\c B = \{ \epsilon_i: = \pi^i x_1 x_2\cdots x_i \}_{i=0}^n$ of the cohomology space $\c W_{\lambda} = H^n(\Omega^{\bullet}_{\c L, \nabla(D)})$. Here we take $\epsilon_0 =1.$ Then $\mathbf{K}_{\lambda} := H^{0}(\Omega^{\bullet}_{\c L^{*},\nabla(D^{*})})$ is a free $\c C_{0}$-module of rank $n+1$ consisting of elements of $\c L^{*}$ which are simultaneously killed by the operators $\{D^{*}_{\ell} \}_{\ell = 0}^n$ and $\mathbf{K}_{\lambda}$ is the $\c C_0$-dual to $\c W_{\lambda}$. It is a free $\c C_0$-module of rank $n+1$ with dual basis $\c B^* = \{ \epsilon^*_i \}_{i=0}^n$, see \cite[section 2]{AS}. 

The connection $\nabla(\lambda \frac{\partial}{\partial \lambda})$ acts  on $\mathbf{K}_{\lambda}$, a free $\c{C}_0$-module of rank $n+1$.  The horizontal sections $\c Y$ satisfy the rank $n+1$ differential system
\begin{equation}
\lambda \frac{\partial}{\partial \lambda}(\c Y(\lambda)) = \c Y(\lambda) G,   
\end{equation}
where $G$ is the  $n+1$ by $n+1$ matrix given by 
\[
G(\lambda) = 
\begin{vmatrix}
0 & 0 & 0 & \dots & 0 & \pi^{n+1}\lambda \\
1 & 0 & 0 & \dots & 0 & 0 \\
0 & 1 & 0 & \dots & 0 & 0 \\                      
\vdots & \vdots & \vdots &\cdots & \vdots & \vdots \\                                           
0 & 0 & 0 & \dots & 1 & 0  \\              
\end{vmatrix}
\ .
\] 
This differential system has a regular singular point with maximal unipotent monodromy at $\lambda = 0$ and an irregular singular point at $\lambda = \infty$. Let $H = G(0)$.
As a consequence  the differential system (1.2) has a fundamental  solution matrix  $\mathbf{Y}_0(\lambda)$ at $\lambda = 0$ consisting of  a full set of log solutions at $\lambda = 0$.  More precisely, a fundamental solution matrix is well-known to have the form 
\begin{equation}
\mathbf{Y}_0(\lambda) = \lambda^H P(\lambda)
\end{equation}
and we may choose the initial condition so that $P(0) =I$. In this situation, the first row vector of $\mathbf{Y}_0(\lambda) = \lambda^H P(\lambda)$ is $P_0(\lambda)$, the first row vector of $P(\lambda)$. It is the unique (up to multiplicative constant) analytic solution at the origin of (1.2) and has the property that at $\lambda = 0$ the initial condition of this row vector is $P_0(0) = (1,0,\ldots,0)$.
There is a Frobenius structure on this differential equation which may be expressed in terms of the Frobenius action on a solution matrix at $\lambda_0$, a unit (or zero) in~$\bar{\Omega}$. Let $\c A$ denote the matrix of Frobenius acting on $\mathbf{K}_{\lambda}$ with respect to the basis ${\c B}^{*}$. The Frobenius structure on a local solution matrix $\mathbf{Y}_{\lambda_0}(\lambda)$ then has the form
\begin{equation}
\mathbf{Y}^{\sigma \phi}_{\lambda_0}(\lambda) \c A(\lambda) = \mathbf{M}_{\lambda_0} \mathbf{Y}_{\lambda_0}(\lambda). 
\end{equation}
where $\mathbf{M}_{\lambda_0}$ is locally constant, $\sigma$ is the Frobenius generator of $\text{Gal}(\Omega/\Omega_1)$ acting on coefficients of series in $\lambda$ and $\phi$ maps the variable $\lambda$ to $\lambda^p$. By a well-known argument, the so-called  ``Dwork trick", we know $\mathbf{Y}_{\lambda_0}(\lambda)$ converges in $D(\lambda_0, 1^{-})$ for any unit $|\lambda_0|=1$. By a similar argument, the same holds for $P(\lambda)$ when $\lambda_0=0$.  At the origin, we have that the first coordinate of the row vector $P_0(\lambda)$, namely $ y = P_{(0,0)}$ satisfies the differential equation
\[
\delta^{n+1}  y = \pi^{n+1}\lambda y
\]
where $ \delta = \lambda \frac{\partial}{\partial \lambda}$.
The recursion from this relation gives the well-known formula
\[
 y(\lambda) = P_{(0,0)}(\lambda) = \sum_{j=0}^{\infty} \frac{\pi^{j(n+1)}}{(j!)^{n+1}} \lambda^j.
\]
Then $P_{(0,0)}(\lambda)$ 
has $p$-integral coefficients as does $P_{(0,i)}(\lambda) = (\delta^{i}  y)(\lambda)$ for $1 \leq i \leq n.$

So the first row of $P(\lambda)$ has $p$-integral coefficients; the first coordinate $P_{(0,0)}(\lambda)$ is a principal unit for all $\lambda \in D(0,1^{-})$; and $P_0(0) = (1,0, \dots, 0) $ The relation (1.4) may be expressed
\begin{equation}
P(\lambda^p) \c A(\lambda) P(\lambda)^{-1} = \lambda^{-pH} \mathbf{M}_0 \lambda^H.          
\end{equation}
It follows easily, since the left side of (1.4) is analytic in $D(0,1^{-})$, that, in fact, no log terms appear on the right. It follows that 
\begin{equation}
P^{\sigma}(\lambda^p) \c A(\lambda) P(\lambda)^{-1} =  \mathbf{M}_0.         
\end{equation}
Finally evaluating at $\lambda = 0$, 
\begin{equation}
\c A(0) = \mathbf{M}_0.                    
\end{equation}
We recall the estimates for the entries $\alpha_{j,i}(\lambda)$ of the matrix $\c A(\lambda)$  when $p>n+1$. Explicitly, for all $\lambda$ with $|\lambda|\leq 1$, 
\begin{align*}
\text{ord}\: \alpha_{j,i}(\lambda)&\geq j \quad\text{for $i,j \in \{0,1,\dots,n\},$}\\
\text{ord} \: \alpha_{j,i}(\lambda) &> j \quad\text{for $i>j$,} \\               
\alpha_{j,j}(\lambda) &= p^j \c U_j(\lambda),
\end{align*}
where $\c U_j (\lambda)$ is a $p$-adic unit, $\c U_j (\lambda) \equiv 1 \  (\text{mod} \ \pi)$. 

 Maximal unipotent monodromy at $\lambda_0 = 0$ enables us to describe $\mathbf{M}_0$ in $(1.4)_{\lambda_0 =0}$ in more detail. Frobenius preserves the subspace consisting of solutions in $D(0,1^{-})$ involving at most the $i$-th power of $\text{log} \ \lambda$. This gives the lower triangularity of $\mathbf{M}_0$. Since $\text{log} \ \lambda^p = p \ \text{log} \ \lambda$, we conclude $\mathbf{M}_{0;j,i} = p\mathbf{M}_{0;j-1,i-1}$. Therefore
\begin{equation}
\mathbf{M}_0 = 
\begin{vmatrix}
1 & 0 & 0 & \dots & 0 & 0 \\a_1  & p & 0 & \dots &   0 & 0 \\ a_2 & pa_1 & p^2 & \dots & 0 & 0 \\  \vdots & \vdots & \vdots & \vdots  & \vdots & \vdots \\                       a_n & pa_{n-1} & p^2a_{n-2} & \dots & \ p^{n-1}a_1 & p^n \\ \end{vmatrix}.
\end{equation}
Here, $a_i = \alpha_{i,0}(0)$. We know from an explicit calculation that $\alpha_{0,0}(0) =1$. Since the unique holomorphic solution at the origin, namely, $\underline{z}(\lambda) = (z_0(\lambda), z_1(\lambda),\dots,z_n(\lambda))$ has $z_0(\lambda) = \displaystyle \sum_{j=0}^{\infty} \frac{\pi^{(n+1)j} \lambda^j}{j!^{n+1}}$ and $z_j(\lambda) = \displaystyle \big(\lambda\frac{\partial}{\partial \lambda}\big)^j(z_0(\lambda))$, it follows that $\underline{z}(0) = (1, 0,\dots,0)$. This is also a consequence of the normalization $P(0) = I$.

 We recall the solutions of the deformation equation (1.2) satisfy a symplectic structure.  
 In the case $p \neq 2$, the choice of $-\pi$ instead of $\pi$ as the chosen $(p-1)$-st root of $-p$ in effect fixes the dual character $\hat{\Theta}$ instead of $\Theta$. 
 We denote wherever possible the effect of making the choice of $-\pi$ by using this as a subscript. Thus, for example, the differential system in $(1.2)$ becomes $(1.2)_{-\pi}$ namely $\lambda\frac{\partial}{\partial \lambda}\c Y(\lambda) =\c Y(\lambda)G_{-\pi}$. 
 Of course, $G_{-\pi} = G$ when $n$ is odd. In the same manner, the Frobenius matrix with respect to the basis $\c B^{*}_{-\pi}$ is denoted $\c A_{-\pi}$, and the particular solution to the differential system (1.2) above at $\lambda_0 = 0$ is given by $\mathbf{Y}_{0, -\pi} = \lambda^H P_{-\pi}(\lambda)$. 
 Then if $\lambda_0$ is a unit and $\mathbf{Y}$ (respectively, $\mathbf{Y}_{-\pi}$) is a  fundamental solution matrix of (1.2) (respectively, $(1.2)_{-\pi}$) then 
\begin{equation}
\mathbf{Y} \Theta \mathbf{Y}_{-\pi}^t = \theta          
\end{equation}
where $\theta$ is constant in $D(\lambda_0, 1^{-})$, and 
\begin{equation}
\Theta = 
\begin{vmatrix}
0 & 0 & 0 & \dots &   0 & 0 & 1 \\
0 & 0 & 0 & \dots &   0 & -1 & 0 \\
0 & 0 & 0 & \dots &   1 & 0 & 0 \\                      
\vdots & \vdots & \vdots & \vdots & \vdots & \vdots & \vdots \\                                    
(-1)^n & 0 & 0 & \dots &  0 & 0 & 0  \\              
\end{vmatrix}.
\end{equation}
Note that the superscript $t$ will denote the matrix transpose throughout this paper. 
At the origin $\lambda_0 =0$, if we take  our solutions to be $ \mathbf{Y}_0 = \lambda^H P(\lambda)$ and $\mathbf{Y}_{0,-\pi} = \lambda^H P_{-\pi}(\lambda)$ in $D(0,1^{-})$ as above, then  $\theta$  satisfies $\theta = \Theta$ and (1.9) becomes in this case 
\begin{equation}
\mathbf{Y}_0 \Theta \mathbf{Y}_{0,-\pi}^t = \Theta.          \end{equation}
We recall as well the functional equation
\begin{equation}
{\c A}  \Theta {\c A^t_{-\pi}}  = p^n\Theta,      
\end{equation}
which holds for $(p, n+1) = 1$. (See Dwork \cite{D3})

The  collection  ($\c A, G(\lambda), \Theta, \underline{\c Y}$) with the properties above will form a  $\textbf{K-Pack\{$N$\}}$ for $N=0$ (which is defined in the next section) and is the starting point of constructing by induction  $\textbf{K-Pack\{$N$\}}$ for $N = \lfloor\frac{n}{2}\rfloor$.  The main applications are then a direct consequence of this result.

\section{Main Result}
In this section we follow Dwork's approach of decomposition by Hodge layers, developed in \cite{D1,D2}. In brief, this decomposition has an arithmetic component which will yield for us $p$-adic analytic formulas for all the eigenvalues of the Frobenius action.  Due to the maximal unipotent monodromy of the deformation equation at the regular singular point $\lambda = 0$, it also has  a classically defined  analytic expression which enables us to characterize these formulas in terms of solutions of this classically defined  confluent hypergeometric differential equation and its exterior powers. We fix $n$ and proceed inductively in terms of an ancillary index $N$ running from $0$ to $\left \lfloor{\frac{n}{2}}\right \rfloor$. We define the required properties of what we call a ``$\textbf{K-Pack\{$N$\}}$'' in I--IV below and then prove their existence via induction on $N$, starting with Lemma 2.1 below.

   
We define a $\textbf{K-Pack\{$N$\}}$  for the family of Kloosterman sums $Kl_{n}$ defined over $\mathbb{G}_m$ to be a collection of $n-2N +1$ by $n-2N +1$ matrices 
\[
\{ \c A(\lambda), G(\lambda),  {\Theta},  \underline{\c Y } \}
\]
as follows.  When needed the matrices in the $\textbf{K-Pack\{$N$\}}$ will have a subscript $N$ added for clarity, so of the form
\[
\{ \c A_N(\lambda), G_N(\lambda),  {\Theta}_N,  \underline{\c Y }_N \}.
\] 
The Frobenius matrix  $\c A$ and the matrix $G$  have entries which are series in $\lambda$ with 
integral coefficients, analytic and bounded by $1$ on $D(0,1^+)$.  The matrix ${\Theta}$ is a constant matrix defined in (\ref{theta}).  We assume that $G(0)$ is maximally nilpotent, $G_N(0) = H_N$ being a square matrix of dimension $n-2N+1$ with ones along the subdiagonal and zeros elsewhere. The associated linear system
\begin{equation}
\lambda \frac{\partial}{\partial \lambda} \c Y(\lambda) = \c Y(\lambda) G_N(\lambda).
\end{equation}
then has maximal unipotent monodromy at the regular singular point $\lambda = 0$. A fundamental solution matrix $\c Y_0$ at $\lambda = 0$ consisting of log solutions may be given by $\c Y_{0} = \lambda^{H_N} P_N(\lambda)$ with $P_N(0) = I_N$. By the so-called  Dwork's trick the entries of $P_N(\lambda)$ converge in the open disk $D(0,1^{-})$. In this definition, we are writing abusively
$\underline {\c Y} $ for a collection $\{ \c Y_{\hat{\lambda}_0} \}_{\hat{\lambda}_0 \in \Gamma}$ of well-behaved (to be specified below in (2.2)) local fundamental solution matrices to (2.1), indexed by the set $\Gamma = \{\text{all Teichm\"uller units}\} \cup \{0\}$, i.e., by $ \Gamma = \{ \hat{\lambda}_0  \in \c O_{\c T} \mid \sigma(\hat{\lambda}_0) = \hat{\lambda}_0^p  \}$. We will write $ \c Y_{\hat{\lambda}_0}$ for the fundamental solution matrix centered at $\hat{\lambda}_0$.  The entries in each matrix $\c Y_{\hat{\lambda}_0}$ for $\hat{\lambda}_0 \neq 0$ are series convergent in $D(\hat{\lambda}_0, 1^{-})$, where $\c T$ is the maximal unramified extension of $\mathbb{Q}_p$ and  $\c U = \c T(\zeta_p)$. We assume that the matrices $\c Y_N$ are normalized in the sense that the row vectors of $\c Y_N$, say $\c Y_N^{(i)}$ with $i$ running from $0 \ \text{to} \  n-2N$ not only are solutions of (2.1) but also are ``eigenvectors" of the semi-linear map from local solutions at $\hat{\lambda}^p_0$ to local solutions at $\hat{\lambda}_0$:
\[
\c Y^{(i)}(\lambda) \rightarrow (\c Y^{(i)})^{\sigma}(\lambda^p)\c A(\lambda)
\]
in the sense that
\begin{equation}
 (\c Y^{(i)}_{\hat{\lambda}_0})^\sigma(\lambda^p) \c A(\lambda) = p^i \c Y^{(i)}_{\hat{\lambda}_0}(\lambda).    
\end{equation}
Finally, if $\rho$ is the element of $\text{Gal} (\mathbb{Q}_q(\pi)/ \mathbb{Q}_q)$ taking $\rho(\pi) = -\pi$, we will denote the image of the action by $\rho$ using the subscript $-\pi$.  The matrices in the {\bf K-Pack{$\{N\}$}} are further required to satisfy Properties I–IV; more precisely, they are required to satisfy (2.3), (2.6), (2.7), (2.8), (2.9),(2.11) and (2.13). (in which we denote the entries of the matrix $\c A(\lambda)$ by $\c A_{j,i}(\lambda))$:

I).
\begin{equation}
\begin{split}
&i). \ \c A(0) \  \text{ is lower triangular with diagonal equal to} \  \c D = \{1, p,\dots,p^{n-2N} \}, \\
   & ii). \  \c A(\lambda) \ \text{factors as follows}: \  \c A(\lambda) =  \c D  \hat{\c A}(\lambda) \ \text{with}\\
    &\hat{\c A}(\lambda) \
    \text{an $n-2N+1$ by $n-2N+1$ matrix with entries }\text{analytic on} \ D(0,1^+) \text{ with integral coefficients. } \\
    &iii). \ \text{and with ord} \ \hat{\c A}_{j,i}(\lambda) > 0 \ \text{for} \ i>j,\\
    &iv). \  \hat{\c A}_{j,j}(\lambda) \ \text{is a principal unit for all} \ j \in \{0, 1,\dots,n-2N\}. 
\end{split}
\end{equation}

It will be useful as well to block subdivide the matrix $\c A$ as follows:
\begin{equation}
\c A= 
\begin{vmatrix}
\c A_{0,0} & \c J_0 & \c A_{0,n-2N} \\
\c T_0 & \c U & \c T_{n-2N}\\
\c A_{n-2N,0} & \c J_{n-2N} & \c A_{n-2N, n-2N} \\      \end{vmatrix}             \hspace{3cm}  
\quad.
\end{equation}
where $\c J_{i}$ denotes the $n-2N-1$ row vector consisting of the elements of row $i$ of the matrix $\c A$ with elements from columns 1 through $n-2N-1$. Similarly, we use $\c T_j$ to denote the elements from $\c A$ from column $j$ and rows $1 $ through $n-2N-1$. Finally, $\c U$ denotes the $n-2N-1$ by $ n-2N-1$ submatrix of $\c A$ having rows $\c J_i$ with $1 \leq i \leq n-2N-1$. We also define a $p$-integral matrix via $\bar{\c U} := p^{-1}\c U$. 

II).\\
There is a Frobenius structure on the differential equation (2.1) which we express in terms of solutions to (2.1) at $\hat{\lambda}_0 \in \Gamma$ as follows.
Let $\hat{\lambda}_0 \in \Gamma$ and $\mathbf{Y}_{\hat{\lambda}_0}(\lambda)$ a fundamental solution matrix of (2.1) at $\hat{\lambda}_0$. Then
\begin{equation}
\mathbf{Y}_{\hat{\lambda}_0}^{\sigma \phi} \c A(\lambda) = \mathbf M_{\hat{\lambda}_0} \mathbf{Y}_{\hat{\lambda}_0}  
\end{equation}
where $\mathbf{M}_{\hat\lambda_0}$ is a constant matrix in a neighborhood of $\hat{\lambda}_0$, $\sigma$ is the Frobenius generator of ${\rm Gal}(\Omega/\Omega_1)$ and where $\phi$ raises the variable $\lambda$ to its $p$-th power and $\sigma$ acts on coefficients. 

From our assumption (2.2),(2.3),(2.5), we deduce that the normalized solution $\c Y$ in the {\bf K-Pack\{$N$\}} at a Teichm\"uller unit $\hat{\lambda}_0$ has the precise form
\begin{equation}
\c Y_{\hat{\lambda}_0}^{\sigma \phi} \c A(\lambda) = \c D \c Y_{\hat{\lambda}_0},
\end{equation}
where $\c D = \c D_N$ is the $n-2N+1$ by $n-2N+1$ diagonal matrix with entries $\{1, p, p^2,\dots,p^{n-2N} \}. $ For the case $\hat{\lambda}_0 = 0$, the Frobenius on the solution matrix $\c Y_{\hat{\lambda}_0}$, has the form
\begin{equation}
\c Y_0^{\sigma \phi} \c A(\lambda) = \c D_0 \c Y_0.   
\end{equation}
where $\c D_0$ is lower triangular and $\c D_0$ and $\c D $ agree on the diagonal. In fact, this follows from our choice of $\c Y_{\hat{\lambda}_0}$ and the Frobenius structure in the case $\hat{\lambda}_0 =0$. We note that as in Section 1, it also follows that $\c A_N(0) = M_{N,0}$ here, where we recall that the $N$ in the subscript on the right refers to the induction index and the $0$ refers to the point in question $\hat{\lambda}_0 =0$. In addition, we make the following requirements (namely, II bis) given below) of the normalized solution $P_N(\lambda)$ at the origin $\hat{\lambda}_0 =0$. These properties will imply the $p$-adic analytic continuation properties of ratios of elements in $P(\lambda)$, see Lemma \ref{P}.
Again, $N$ is the induction index. More precisely, we require

II bis).

\begin{equation}
\begin{split}
&i).\ P_N(0) =I. \\
&ii).\ P_N(\lambda)  \ \text{converges in the open disk} \ D(0,1^{-}).
 \end{split}
\end{equation}
The same requirements are imposed on the system obtained by replacing
\(\pi\) by \(-\pi\).

III). \\
We assume as well a symplectic relation satisfied by solutions, say, at $\hat{\lambda}_0 \in \Gamma$
\begin{equation}
\mathbf{Y} \Theta \mathbf{Y}_{-\pi}^t = \theta.
\end{equation}
Here 
\begin{equation}
\Theta_N =
\begin{vmatrix}\label{theta}
0 & 0 & 0 & \cdots &\ 0 & 1\\
0 & 0 & 0 & \cdots &-1 & 0\\
\vdots & \vdots & \vdots & \cdots & \vdots & \vdots \\       0 & (-1)^{n-2N-1} & 0 & \cdots & 0 & 0\\
(-1)^{n-2N} & 0 & 0 & \cdots & 0 &  0 
\end{vmatrix}                      
\end{equation}
 and $\theta$ is a constant matrix. In the special case where $\hat{\lambda}_0 = 0$ and where the fundamental solution matrix is given by (1.3) we assume
\begin{equation}
\mathbf{Y}_0 \Theta_N \mathbf{Y}_{0,-\pi}^t = \Theta_N.   
\end{equation}
We will denote the $n-2N-1$ by $n-2N-1$ matrix  consisting of the entries of $\Theta_N$ in the rows and columns  $1$ through $n-2N-1$ by $\Theta^{(1)}_N$, so that 
$(\Theta^{(1)}_N)_{ij} = (\Theta_N)_{i+1,j+1}$ with $0 \leq i,j \leq  n-2N-2$.   
For convenience, we set $\epsilon = \epsilon_N = (-1)^{n-2N}$. Then the matrix $\Theta_N$ has the block form
\begin{equation}
\begin{vmatrix}
0 & 0  & 1\\
0 & \Theta^{(1)}_N & 0 \\
\epsilon_N & 0 & 0  
\end{vmatrix} .
\end{equation}
 Clearly $\Theta_{N+1} = - \Theta^{(1)}_{N} $.
 
IV). \\
Finally we assume a functional equation of the form
\begin{equation}
\c A(\lambda) \Theta \c A_{-\pi}^t(\lambda) = p^{n-2N} \Theta.
\end{equation}

\begin{lemma}\label{P}
Let
\[
Z(\lambda)
=
\bigl(Z_0(\lambda),\ldots,Z_{n-2N}(\lambda)\bigr)
\]
be a vector of analytic function on $D(0,1^{-})$. Assume the following conditions hold:
\begin{enumerate}
\item[(1)] $Z(0)=(1,0,\cdots,0)$;
\item[(2)] $Z^{\sigma}(\lambda^p)\c A(\lambda)=Z(\lambda)$, where $\sigma$ is the Frobenius generator acting on coefficients, $\c A$ is the Frobenius matrix defined in a $\textbf{K-Pack\{$N$\}}$.
\end{enumerate}
Then,\[
Z_{0}(\lambda)\in 1+\pi\mathcal O_{\mathcal U}[[\lambda]],
\qquad
Z_{j}(\lambda)\in\pi\mathcal O_{\mathcal U}[[\lambda]]
\quad\text{for }j>0.
\]
\end{lemma}
\begin{proof}
Write
\[
Z_{j}(\lambda)=\sum_{m\geq0}p_{j,m}\lambda^m,
\qquad
A(\lambda)=\bigl(A_{ij}(\lambda)\bigr)_{0\leq i,j\leq n-2N},
\qquad
A_{ij}(\lambda)=\sum_{b\geq0}a_{ij;b}\lambda^b.
\]
The $j$-th coordinate of the equation in (2) is
\[
\sum_{i=0}^{n-2N}Z_{i}^{\sigma}(\lambda^p)A_{ij}(\lambda)
=
Z_{j}(\lambda).
\]
Comparing the coefficients of $\lambda^m$, we obtain
\[
p_{j,m}
=
\sum_{i=0}^{n-2N}
\ \sum_{pa+b=m}
\sigma(p_{i,a})a_{ij;b}.
\]
Since $p_{j,0}\in\mathcal O_{\mathcal U}$ for every $j$ and
$a_{ij;b}\in\mathcal O_{\mathcal U}$, the right-hand side for $m>0$
involves only the coefficients $p_{i,a}$ with $a<m$. It follows by
induction on $m$ that
$p_{j,m}\in\mathcal O_{\mathcal U}$ for all $j$ and $m\geq0.$ By Condition (2.3), $\c A(\lambda)\equiv E_{00} \text{ mod } \pi$, where $E_{00}$ denote the matrix whose \((0,0)\)-entry is \(1\) and all other entries are \(0\).
Reducing the equation modulo $\pi$, we obtain
\[
Z^{\sigma}(\lambda^p)E_{00}
\equiv
Z(\lambda)
\pmod{\pi}.
\]
It follows immediately that $Z_j(\lambda)\equiv0\pmod{\pi}$
for $j>0$, and
$Z_0^{\sigma}(\lambda^p)
\equiv
Z_0(\lambda)
\pmod{\pi}.$ Suppose that $m>0$ is the smallest integer such that
$p_{0,m}\not\equiv0\pmod{\pi}$. The last congruence implies successively
that $p^k\mid m$ for every $k>0$, which is impossible. Hence
$Z_0(\lambda)$ is congruent modulo $\pi$ to a nonzero constant which is $1$ by our assumption \rm{(1)}. The result follows.
\end{proof}

We begin our induction by showing the existence of a $\textbf{K-Pack\{$0$\}}$.
In the case of the $n$-dimensional hyperkloos\-terman family (with $p>n+1$) of  the introduction, with deformation equation (1.2), we choose normalized solution matrices to (1.2) by specifying the rows $\c Y^{(i)}_{\hat{\lambda}_0}$ of the  solution matrix $\c Y_{\hat{\lambda}_0}$ as in (2.2) above, for $\hat{\lambda}_0 \in \Gamma^{*}, (\Gamma^{*} = \Gamma - \{0\})$ and $i=0,\dots,n$. Similarly we take $\c Y_{\hat{\lambda}_0}(\lambda) = \lambda^H P(\lambda)$ at the origin.

\begin{lemma}
Using this choice of normalized solution matrices $\c Y$ and the matrices $\c A, G, {\Theta}$ from the preceding section, we obtain a {\rm \textbf{K-Pack\{$0$\}}}.
\end{lemma}

\begin{proof}
At a unit, $\hat{\lambda}_0 \in \Gamma$,  we know from the estimates of Frobenius in Section 1 and the  work of Dieudonn\'e and Manin, that there exist vectors $\tilde{\xi}_i \in \c O_{\c U}^{n+1}$ which are ``eigenvectors" under the the semilinear map
\[
\xi \rightarrow \sigma(\xi) \c A(\hat{\lambda}_0)
\]
(with $i=0, \dots, n$) satisfying 
\[
\tilde{\xi}_i^{\sigma} \c A(\hat{\lambda}_0) = p^i \tilde{\xi}_i.
\]
It follows then from the commutativity of Frobenius and Dwork's deformation map that there is a series, say $\xi_i(\lambda)$  in $\c U[[\lambda - \hat{\lambda}_0]]$, with $\xi_i(\hat{\lambda}_0) = \tilde{\xi}_i$ and
\[
\xi_i^{\sigma}(\lambda^p) \c A(\lambda) =p^i \xi_i(\lambda)
\]
in $D(\hat{\lambda}_0, 1^{-}).$ We define the normalized solution matrix $\c Y_{\hat{\lambda}_0}$ to be the matrix with rows $\{ \xi_0, \xi_1,\dots,\xi_{n} \}$. Similarly, at the  point $0$, we take $\c Y_0 = \lambda^H P(\lambda) $ with $P(0) =I$. Note that in the open unit neighborhoods of $\hat{\lambda}_0 = 0$ and $\hat{\lambda}_0$ a Teichm\"uller unit, the first row of each solution matrix is a fixed point of the Frobenius map 
\[
Y(\lambda) \rightarrow Y^{\sigma \phi} \c A(\lambda).
\]
This vector has all entries bounded by 1 in the open unit disk at $\hat{\lambda}_0$ and the first entry is a unit there. (Of course in the case $N=0$ which we are treating here, the results of $ \text{II}\ \text{bis})_{N=0}$ are immediate from the well-known series formulas for $\c Y_0(\lambda)$.) The other required properties to establish that this gives a {\bf K-Pack\{$0$\}} are from the  references in the first section.
\end{proof}

\begin{lemma} 
Given a {\rm \textbf{K-Pack\{$N$\}}}, $0 \leq N \le \left \lfloor{\frac{n}{2}}\right \rfloor$ with
\[
\{ \c A_N(\lambda), G_N(\lambda),  {\Theta}_N,  \c Y_N  \},
\] 
then for any Teichm\"uller unit $\hat{\lambda}_0 \in \Gamma$ with $\sigma^a(\hat{\lambda}_0) = \hat{\lambda}_0$ the matrix $\c A_N(\hat{\lambda}_0^{\frac{q}{p}})\cdots\c A_N(\hat{\lambda}_0)$ has a unique unit eigenvalue, say $\omega_N^{(0)}({\lambda}_0)$, where $\lambda_0$ is the reduction mod $\pi$ of $\hat{\lambda}_0$. It also has, by the functional equation $(2.11)$, a unique eigenvalue 
\begin{equation}
\omega_N^{(n-2N)}({\lambda}_0) = \frac{q^{n-2N}}{\omega_{N,-\pi}^{(0)}({\lambda}_0)} 
\end{equation}
of maximal ${\rm ord}_q = n-2N$. Let $\c Y^{(N)}_{ \hat{\lambda}_0}$ be chosen at $\hat{\lambda}_0 = 0$ to be $\c Y^{(N)}_ {\hat{\lambda}_0} = \lambda^{H_N}P_N(\lambda)$ with $P_N(0) = I_{n-2N+1}$. For notational convenience let $\c F(\lambda) = P_{0,0}(\lambda)$ where we have omitted the induction index  $N$. Substituting the first row of \(P(\lambda)\) for \(Z(\lambda)\) in Lemma 2.1, $\c F(\lambda)$ is a principal unit for all $\lambda\in D(0,1^{-})$, The ratios $\eta_{0,j}(\lambda) = \frac{P_{0,j}(\lambda)}{P_{0,0}(\lambda)} \ \text{and  the ratio} \ \frac{\c F(\lambda)}{\c F(\lambda^q)}$, which are well-defined by Lemma \ref{P} and integral in this disk, in addition have analytic continuation to the disk $D(0,1^{+})$.
The following p-adic analytic formulas hold then for the distinguished eigenvalues
\[
\omega_N^{(0)}({\lambda}_0) = \frac{\c F_N(\lambda)}{\c F_N(\lambda^q)}|_{\lambda = \hat{\lambda}_0}
\]
and a similar formula (using $(2.11)$) for the unique eigenvalue of maximal ${\rm ord}_q = n-2N$, namely $\omega_N^{(n-2N)}({\lambda}_0) = \frac{q^{n-2N}}{\omega_{N,-\pi}^{(0)}({\lambda}_0)}$. 
\end{lemma}

Before proving this lemma, we mention a useful immediate consequence.
\begin{corollary} 
Set $\underline{\eta}_1 = (\eta_{0,1},\dots,\eta_{0,n-2N-1})$ and $\underline{\eta}_2 = \eta_{0, n-2N}.$  In addition, set $\underline{\eta} =  (\underline{\eta}_1, \underline{\eta}_2)$.  Define the  matrix $\c Y_2(\lambda)$  in block matrix notation as follows:
\begin{equation}
\c Y_2 =
\begin{vmatrix}
1 & \underline{\eta}_{1} & \underline\eta_{2}\\
0 & I & -\epsilon {\Theta}^{(1)}\underline{\eta}_{1,-\pi}^t\\                                                     0 & 0 & 1
\end{vmatrix},                 
\end{equation}
where the matrix has been block divided with rows and columns having dimensions $\{1, n-2N-1, 1\}$. Then the entries converge in $D(0,1^+)$ and are bounded by $1$ there.  The matrix $\c Y_2$  is invertible in the ring of $n-2N+1$ by $n-2N+1$ matrix functions defined and analytic on $D(0,1^+)$ and bounded by $1$ there.
\end{corollary}
The corollary follows directly from the lemma. We next prove the lemma.

\begin{proof}
We note first of all that proceeding at $\hat{\lambda}_0 =0$ from the relationship   
\begin{equation}
(\c Y_N^{(0)})^{\sigma \phi} \c A_N = \c Y_N^{(0)}           \end{equation}
and using from II bis) the initial data $\c Y_N^{(0)}(0) = (1,0,\dots,0)$, the coefficients of the series $\c Y_N$ in $\lambda$ are  $p$-integral and  $\c Y^{(0)}_{N; 0,0}(\lambda) $ is a principal unit for $\lambda \in D(0,1^{-})$. So the vector converges in the open disk $D(0, 1^{-})$ and $\c F_N = \c Y^{(0)}_{N;0,0}$ is a unit there. We use a fixed point argument due to Dwork.  Given the properties of the matrix $\c A_N$ assumed in (2.3) it is a well-known argument that the map $\mathbf{\Phi}$ defined in (2.17) below which takes $(V_1(\lambda), V_2(\lambda),\dots,V_{n-2N}(\lambda))$ to
 \begin{equation}
 \bigg(\frac{\c A_{0,1} + V_1^\sigma(\lambda^p) \c A_{1,1} +\cdots+V_{n-2N}^\sigma(\lambda^p) \c A_{n-2N,1}}{\c A_{0,0} + V_1^\sigma(\lambda^p) \c A_{1,0} +\cdots+V_{n-2N}^\sigma(\lambda^p) \c A_{n-2N,0}},\dots,\frac{\c A_{0,n-2N} + V_1^\sigma(\lambda^p) \c A_{1,n-2N} +\cdots+V^\sigma_{n-2N}(\lambda^p) \c A_{n-2N,n-2N}}{\c A_{0,0} + V_1^\sigma(\lambda^p) \c A_{1,0} +\cdots+V_{n-2N}^\sigma(\lambda^p) \c A_{n-2N,0}}\bigg)
 \end{equation}  
is a contraction map on the complete metric space $\c S$ consisting of $(n-2N)$-tuples $(V_1(\lambda), \dots, V_{n-2N}(\lambda))$ of analytic functions convergent in $D(0,1^{-})$ and having sup norm on $D(0,1^{-}) \leq  |\pi|$. Since $(\c Y_N^{(0)})^{\sigma \phi} \c A_N = \c Y_N^{(0)}$, it follows then, dividing both sides by $\c Y_{N;0,0}$, that $\underline{\eta}$ is this fixed point. 
In the same way, $\mathbf{\Phi}$ is  a contraction map as well on the complete metric space $\hat{\c S}$ consisting of $(n-2N)$-tuples $(V_1, \dots, V_{n-2N})$ of elements defined and analytic on $D(0,1^{+})$ and having sup norm $\leq |\pi|$. In particular, this proves that the unique fixed point $\underline{\eta}$  lies in $\hat{\c S}$. It follows from (2.16) that
\begin{equation}
(1, \underline{\eta})^{\sigma \phi} \c A = \frac{\c F_N(\lambda)}{\c F_N^\sigma(\lambda^p)}(1, \underline{\eta}),
\end{equation}
from which the assertion concerning the unit eigenvalue follows.   The assertion concerning the eigenvalue of maximal ord then follows from the functional equation (2.11).
\end{proof}

\begin{theorem}
Let $\hat{\lambda}_0 \in \Gamma$ with reduction $\text{mod} \ \pi = \lambda_0$ satisfy $\sigma^a(\hat\lambda_0)=\hat\lambda_0$ . For any $N$ we denote the set of eigenvalues of $\c A_N(\hat{\lambda}_0^{\frac{q}{p}})\cdots\c A_N(\hat{\lambda}_0)$ by
\[
\c E_N(\lambda_0) = \{ \omega_N^{(i)}(\lambda_0) \}_{i=0}^{n-2N}.
\]
Given a {\rm \textbf{K-Pack\{$N$\}}} for the sums $Kl_{n}$,
\[
\{ \c A_N(\lambda), G_N(\lambda),  \mathbf{\Theta}_N,  \c Y_N  \}
\]
there is a {\rm \textbf{K-Pack\{$N+1$\}}}  such that the relationship between $\c E_N(\lambda_0) \ \text{and} \ \c E_{N+1}(\lambda_0) $ is given by
\[
\c E_N(\lambda_0) = q \c E_{N+1}(\lambda_0) \cup \{\omega_N^{(0)}(\lambda_0), \omega_N^{(n-2N)}(\lambda_0) \},
\] 
where the first term on the right signifies multiplying every element of the set $\c E_{N+1}$ by $q$.
\end{theorem}

\begin{proof}
To simplify the notation $\rm \textbf{K-Pack\{$N+1$\}}= \{ \c A_{N+1}, G_{N+1}, \Theta_{N+1}, \c Y_{N+1} \}$ when $N$ is understood, we will often write $\rm \textbf{K-Pack\{$N+1$\}}= \{ \c A^{(1)}, G^{(1)}, \Theta^{(1)}, \c Y^{(1)} \}$. 
Much of this follows along the same lines as Dwork's proof of Theorem 1.1 in \cite{D2}. To the extent that some parts of his argument are repeated here, the purpose is to make the exposition more self-contained. 

We use the ``global" factor $\c Y_2$ defined in the corollary above to factor the local matrices in the collection $\c Y_N$. This gives a collection of local matrices $\c Y_1 = \c Y \c Y_2^{-1}$ where $\c Y_1$ is notational shorthand for the collection $\{ \c Y_{N; 1; \hat{\lambda}_0} \}_{\hat{\lambda}_0 \in \Gamma}$.  We assert that each matrix $\c Y_{1; \hat{\lambda}_0}$ is block lower triangular and may be written in the form
\begin{equation}
\c Y_1 =
\begin{vmatrix}
U & 0 & 0\\
\underline\zeta U & \c Y^{(1)} & 0\\     TU & \underline z_2 & z_3\\
\end{vmatrix}                 
\end{equation}
All of the functions appearing in the matrix $\c Y_1$ are defined below in Lemma \ref{z23}.  We have for the {\bf K-Pack\{$N$\}}  the following relations.  We note (as observed in the previous section) that indeed they do hold for $N=0$. We use the asterisk superscript as in $Z^{*}$ applied to a matrix $Z$ to denote $(Z^{-1})^t$. Then the following relations hold.
\begin{align}
&{\c Y} {\Theta} = \theta {\c Y}^{*}_{-\pi},    \\
&{\c Y}^{\sigma \phi}_{\hat{\lambda}_0} {\c A} {\c Y}^{-1}_{\hat{\lambda}_0} = {\c D} \  \text{for all units}  \ \hat{\lambda}_0  \in \Gamma,\\
&{\c Y}_0^{\sigma \phi} {\c A} {\c Y}_0^{-1} = {\c D}_0,   \\   
&{\c A} {\Theta} {\c A}_{-\pi}^t = p^{n-2N} {\Theta}.   
\end{align}

Relation (2.21)  also has a dual version
\[
\c Y^{\sigma \phi}_{-\pi, \hat{\lambda}_0} \c A_{-\pi} \c Y^{-1}_{-\pi,\hat{\lambda}_0} = {\c D} \ \text{for all} \ \hat{\lambda}_0 \in \Gamma\setminus\{0\},
\]
(and (2.22), for $\hat{\lambda}_0 =0$ as well).
It follows then from relation (2.21) that 
\[\c D \theta \c D^t= \c Y^{\sigma \phi} \c A \c Y^{-1} \theta (\c Y_{-\pi}^{-1})^t \c A^t_{-\pi} (\c Y_{-\pi}^{\sigma \phi})^t,
\]  
from which using relation (2.20) we get the equality
\[
\c D \theta \c D^t = \c Y^{\sigma \phi} \c A {\Theta} \c A_{-\pi}^t  (\c Y_{-\pi}^{\sigma \phi})^t
\]
from which we get by (2.23)
\[\c D \theta \c D^t = p^{n-2N} \c Y^{\sigma \phi}  \Theta (\c Y_{-\pi}^{\sigma \phi})^t
\]
and finally using (2.20) again
\begin{equation}
\c D \theta \c D^t= p^{n-2N} \theta^{\sigma}    
\end{equation}
follows immediately. The relation (2.22) implies that the matrix $\theta$ is antidiagonal with entries $\theta_{ij} = 0$ whenever $i+j \neq N-1$ and $\theta_{i,N-i-1}^{\sigma} = \theta_{i, N-i-1}$ for $i = 0,1,\dots,N-1$. Let us denote for any unit $\hat{\lambda}_0 \in \Gamma$ the first row of our normalized solution by $(U_{0,0}, \underline{U}_{0,1}, U_{0,2}) = U_{0,0}(1,\underline{\eta}_1, \eta_2)$. It follows that $U_{0,0}$ is a unit in the open disk centered at $\hat{\lambda}_0$. Since $\theta$ is antidiagonal, it follows from (2.20) above that
\[
(U_{0,0}, \underline{U}_{0,1}, U_{0,2}){\Theta}(U_{-\pi,0,0}, \underline{U}_{-\pi,0,1}, U_{-\pi, 0,2})^t = 0.
\]
From this it follows that 
\begin{equation}
\epsilon \eta_{2} + \underline{\eta}_{1} {\Theta}^{(1)} \underline{\eta}_{-\pi, 1}^t + \eta_{-\pi, 2} = 0.   
\end{equation}
It follows from  a direct calculation that $\c Y_2 {\Theta} \c Y_{2, -\pi}^t$ agrees with ${\Theta}$ except possibly in the (0,0) coordinate. But this term  is given by (2.25) so that
\begin{equation}
\c Y_2 {\Theta} \c Y_{2, -\pi}^t = {\Theta}.                 \end{equation}
By definition,  
\begin{equation}
\c Y_{1,\hat{\lambda}_0} = \c Y_{\hat{\lambda}_0} \c Y_2^{-1}.                    
\end{equation}
But then (2.26) and (2.20) imply
\begin{equation}
\c Y_1 \mathbf{\Theta} \c Y_{1, -\pi}^t = \theta.           \end{equation}
Let us write $U = U_{0,0}$ when convenient. We have the following lemma:
\begin{lemma}\label{z23}
$\c Y_1$ has the following form which is block lower triangular.
\begin{equation} 
\c Y_1 =
\begin{vmatrix}
U & 0 & 0\\
\underline\zeta U & \c Y^{(1)} & 0\\   TU & \underline z_2 &z_3\\
\end{vmatrix}.            
\end{equation}
Here $z_3=\epsilon\theta_{n-2N,0}U_{-\pi}^{-1},$ $\underline{z}_2 = -\theta_{n-2N,0} \underline{\zeta}^t_{-\pi}({\b\Theta}^{(1)} ({\c Y}^{(1)}_{-\pi})^t)^{-1}$,  where $\theta_{n-2N,0}$ is the $(n-2N,0)$-entry of $\theta$ .
\end{lemma}
\begin{proof}
Let $e_0=(1,0,\cdots,0)$. Define 
$$\bar{A} = {\c Y}_2^{\sigma \phi}\c {A}{\c Y}_2^{-1}.$$
    By (2,18), we have
$$
e_0\bar{A}=e_0{\c Y}_2^{\sigma \phi}\c {A}{\c Y}_2^{-1}
=(1,\underline \eta)^{\sigma\phi}\c {A}{\c Y}_2^{-1}
=\frac{\c F(\lambda)}{\c F^\sigma(\lambda^p)}(1,\underline\eta){\c Y}_2^{-1}
=\frac{\c F(\lambda)}{\c F^\sigma(\lambda^p)} e_0.
$$
Thus, $\bar{A}$ has the form 
$$\begin{vmatrix} a& 0\\ C&B\\ \end{vmatrix},$$
where $a=\frac{\c F(\lambda)}{\c F^\sigma(\lambda^p)}$, $C,B$ are block submatrices. From $\c A(\lambda)\equiv E_{00} \text{ (mod } \pi)$ we deduce that $\bar{A}(\lambda)\equiv E_{00} \text{ (mod } \pi)$ and consequently $B\equiv 0 \text{ (mod } \pi)$. From the Frobenius equation (2.6), we deduce that
$$\c Y_{1,\hat{\lambda}_0}^{\sigma \phi} \bar {A}(\lambda) = \c D \c Y_{1, \hat{\lambda}_0}.$$
Write the first line of $\c Y_{1,\hat{\lambda}_0}$ by $(U,W)$, we have
$$W^{\sigma\phi}B=W.$$
Taking the Gauss norm gives $|W|=|W^{\sigma\phi}B|\le |\pi||W|$, which means $W=0$.
Since $U$ is invertible, we write $\c Y_1$ by
$$\c Y_1 =
\begin{vmatrix}
U & 0 & 0\\
\underline\zeta U & \c Y^{(1)} & \underline z_1\\   TU & \underline z_2 & z_3\\
\end{vmatrix}.     $$
The $(0,1)$-block of (2.27) gives
$$U\underline z_{1,-\pi}^t=0.$$
We get $\underline z_1=0$. The $(n-2N,0)$-entry gives
$$z_3=\epsilon\theta_{n-2N,0}U_{-\pi}^{-1}.$$
Finally, comparison of the \((n-2N,1)\)-block gives
$$\underline{z}_2 = -\theta_{n-2N,0} \underline{\zeta}^t_{-\pi}({\Theta}^{(1)} ({\c Y}^{(1)}_{-\pi})^t)^{-1}.$$
\end{proof}

Looking ahead, we will show $\c Y^{(1)}$ is a family of period matrices for a differential equation of rank $n-2N -1$ and we will build the \textbf{K-Pack}$\{N+1\}$ with $\c Y_{N+1} := \c Y^{(1)}$. 

We turn now to the Frobenius for the \textbf{K-Pack}$\{N+1\}$. From the equation
\[
\bar{A}=(\c Y_1^{\sigma\phi})^{-1} \c D \c Y_1.
\]
We obtain $\bar{A}$ is lower block triangular and may be written explicitly  for a column vector $\bar{\c T}_0$
and row vector $\bar{\c J}_{n-2N}$ both of length $n-2N-1$ and a square matrix $\bar{\c V}$ of dimension $n-2N-1$ as follows:
\begin{equation}
\bar{A} =
\begin{vmatrix}
\bar{A}_{0,0} & 0 & 0\\
\bar{\c T}_{0} & \bar{\c V} & 0\\   \bar{A}_{n-2N,0} & \bar{\c J}_{n-2N} & \bar{A}_{n-2N,n-2N}\\
\end{vmatrix} .                    
\end{equation}
We will give the block entries of (2.30) more explicitly using (2.4) and (2.29). Let ${\c D}^{(1)}$ to be  the $n-2N-1$ by $n-2N-1$ diagonal matrix with entries $\{ 1, p,\dots,p^{n-2N-2} \}$.
Clearly from the estimates in (2.3) and (2.4), and using Lemma \ref{P}, there are  matrices $\hat{\c T}_{0},  \hat{\c J}_{n-2N},   \hat{\c V}, \hat{A}_{n-2N,0}, \hat{A}_{n-2N, n-2N}, \hat{A}_{0.0} \ \text{and} \ \c A^{(1)}$ with integral coefficients, analytic and bounded by $1$ on the closed unit disk $|\lambda|\leq1$,  such that 
\begin{equation}
\bar{A} =
\begin{vmatrix}
\bar{A}_{0,0} & 0 & 0\\
p{\c D}^{(1)}\hat{\c T}_{0} & p  \c A^{(1)} & 0\\           p^{n-2N}\hat{A}_{n-2N,0} & p^{n-2N}\hat{\c J}_{n-2N} & p^{n-2N} (\bar{A}_{(0,0),-\pi})^{-1}\\
\end{vmatrix}.                  
\end{equation}
Moreover, that there is a  square $n-2N-1$ by $n-2N-1$ matrix with entries analytic on $D(0,1^+)$ and bounded by $1$, say, $\tilde{\c A}^{(1)}$,  such that $\c A^{(1)} = \c D^{(1)} \tilde{\c A}^{(1)}$. 

We record a few useful consequences of the above construction.  We note that 
\begin{equation}
(U^{\sigma \phi})^{-1}U = \bar{A}_{0,0},    
\end{equation}
and 
\begin{equation}
({\c Y}^{(1})^{\sigma \phi} \c A^{(1)} = {\c D}^{(1)} {\c Y}^{(1)}.  
\end{equation}
The matrix $\c A^{(1)}$ is the Frobenius in the \textbf{K-Pack}$\{N+1\}$. 
It is worthwhile to get more precision about the ``middle" block matrix $\c Y^{(1)}$, particularly at the origin. We can calculate this matrix, say at $\hat \lambda_0 = 0$, by proceeding from (2.25), namely multiplying the middle block of rows of $\c Y = \lambda^H P(\lambda)$ by the middle block of columns of $\c Y_2^{-1}$. The middle block of rows of $\c Y_{\hat{\lambda} = 0}$ is the product of the middle block of rows of $\lambda^H$ times the matrix $P(\lambda)$. We have first of all
\begin{equation}
\lambda^H =
\begin{vmatrix}
1 & 0 & 0\\
\bar{\c N} & \lambda^{H^{(1)}} & 0\\                        \frac{(\log\lambda)^{n-2N}}{(n-2N)!} & \underline{\c N} & 1\\
\end{vmatrix},                  
\end{equation} 
where we have denoted $\bar{\c N}$ to be the column vector with entries $\{ \log \lambda, \frac{\log^2 \lambda}{2!}, \dots,
\frac{\log^{n-2N-1} \lambda}{(n-2N-1)!} \}$  in that order. Similarly, $\underline{\c N} $ is the row vector with entries $\{\frac{\log^{n-2N-1} \lambda}{(n-2N-1)!},   \dots, \frac{\log^2 \lambda}{2!},  \log\lambda \}$ in that order. The matrix $H^{(1)}$ is the $n-2N-1$ by $n-2N-1$ ``shift" operator.
It is useful to write the block decomposition of the matrix $P(\lambda)$ explicitly as follows: 
\begin{equation}
P =
\begin{vmatrix}
P_{0,0} & \underline{P}_{0} & P_{0,n-2N}\\
\bar{P}_0 &  P^{(1)} & \bar{P}_{n-2N}\\                   P_{n-2N,0} & \underline{P}_{n-2N} & P_{n-2N,n-2N}\\
\end{vmatrix}.                  
\end{equation}
Here, for $0 \leq i\le n-2N$, we have used the notation $\bar{P}_i$ to denote the column vector $\bar{P}_i=(P_{1,i}, \dots, P_{n-2N-2, i})^t$ and similarly the row vector $\underline{P}_i =(P_{i,1}, \dots, P_{i ,n-2N-2})$. Finally, the matrix $P^{(1)}=(P_{i,j})$ with indices $i,j= 1, \dots, n-2N-1$.
We denote the middle block-column of $\c Y_2^{-1}$  by  $\bar{\c Q}= (-\underline{\eta}_1, I,0)^t$. It follows then that
\begin{equation}
\c Y^{(1)} = (\bar{\c N}, \lambda^{H^{(1)}}, 0)  P  \bar{\c Q} = \bar{\c N}(- P_{0,0} \underline{\eta}_1 + \underline{P}_0)  + \lambda^{H^{(1)}}(-\bar{P}_0\underline{\eta}_1 +  P^{(1)}) = \lambda^{H^{(1)}} \c P^{(1)}(\lambda),
\end{equation}
where we have taken 
$\c P^{(1)} = (-\bar{P}_0\underline{\eta}_1 +  P^{(1)})$ and we use that the coefficient of $\bar{\c N}$ in (2.36) above, namely $(- P_{0,0} \underline{\eta}_1 + \underline{P}_0) $ is zero, from the definition of $\underline{\eta}_1$. Hence $\c P^{(1)}(0) = I$ since $P(0) =I$. We take $P_{N+1} = \c P^{(1)}$. It follows that at the origin, $\hat{\lambda} = 0, \ P_{N+1}(0) =  P^{(1)}(0) = I $, i.\ e.,  property (i) of II bis) holds for index $N+1$. 

Substituting the expression in (2.36) into (2.33) and bringing all log terms to one side, we obtain
\begin{equation}
(\c P^{(1)})^\sigma(\lambda^p) \c A^{(1)} (\c P^{(1)}(\lambda))^{-1} = \lambda^{-pH}\c A^{(1)}(0)\lambda^H.
\end{equation}
Then the left side is uniform near the origin so the log terms on the right cancel and we get
\begin{equation}
  (\c P^{(1)})^\sigma(\lambda^p) \c A^{(1)} = \c A^{(1)}(0) \c P^{(1)}(\lambda).
\end{equation}
The convergence of $\c P^{(1)}=P^{(1)}-\bar P_0\underline\eta_1$ in $D(0,1^{-})$ follows directly from that of $P$ and $\underline\eta_1$. Thus property (ii) of II bis) holds for index $N+1$; the preceding identity gives the Frobenius relation at the origin.

We also wish to use this calculation to obtain information on the first row of $\c Y^{(1)}$. We show below that $\c Y^{(1)}$ is a fundamental solution matrix of the equation (2.45) below which arises in the \textbf{K-Pack}$\{N+1\}$. In fact, by (2.33), it is a normalized period solution matrix with the $i$-th row $\c Y^{(1)}_i$ a (semi-linear) eigenvector of $\c A^{(1)}$ corresponding to the (semilinear) eigenvalue $p^i$.  Then using the subscript notation $P_{N+1}:= \c P^{(1)}$ and  (2.36), we conclude that 
the first row of $\c Y^{(1)}$, say $\c Y^{(1)}_0$, is in fact the first row of $\c P^{(1)}$, say $\c P^{(1)}_0$, which is
\[ 
\c P^{(1)}_0 = (-P_{1,0}  \underline{\eta}_1 +  \text{ the  first row  of }   P^{(1)}). 
\]

The relations between the entries of $\bar{A} $ (given in (2.31) above) and $\c {A}$ which follow from (2.29) are listed here
\begin{equation}
\begin{split}
i).&\ \bar{A}_{0,0} = \c A_{0,0} +  \underline{\eta}_{1}^{\sigma \phi} \c T_{0} + \underline{\eta}_{2}^{\sigma \phi}\c A_{n-2N,0}  \\
ii).&\ \bar{T}_{0} = \c T_{0} - \epsilon {\bold\Theta}^{(1)}(\underline{\eta}_{1, -\pi}^{\sigma \phi})^t \c A_{n-2N,0}  \\
iii).&\ \bar{A}_{n-2N,0} =\c A_{n-2N,0}  \\
iv).&\ \bar{\c V}= p\c A^{(1)} = \c U -\c T_{0}\underline{\eta}_{1} + \epsilon  {\bold\Theta}^{(1)}(\underline{\eta}_{1, -\pi}^{\sigma \phi})^t(\c A_{n-2N,0}\underline{\eta}_{1} - \c J_{n-2N}) \\
v).&\ \bar{\c J}_{n-2N  } = \c J_{n-2N} - \c  A_{n-2N,0}\underline{\eta}_{1}.           
\end{split}
\end{equation}

From (2.3) and (2.4) there are $p$-integral matrices $\hat{U}$ and $\tilde{U}$ such that $U = p \tilde{U} = p \c D^{(1)} \hat{U}$. If we use this and the estimates above then the relations above take the following form with $p$-integral matrices on the right side.
\begin{equation}
\begin{split}
i').&\ \bar{A}_{0,0} = \c A_{0,0} + p \underline{\eta}_{1}^{\sigma \phi} {\c D}^{(1)}\hat{\c T}_{0} + p^{n-2N} \underline{\eta}_{2}^{\sigma \phi}\hat{A}_{n-2N,0}  \\
ii').&\ \bar{T}_{0} = p\c 
D^{(1)} \hat{T}_{0} - p^{n-2N}\epsilon {\bold\Theta}^{(1)}(\underline{\eta}_{1, -\pi}^{\sigma \phi})^t \hat{A}_{n-2N,0}  \\
iii').&\ \bar{A}_{n-2N,0} = p^{n-2N} \hat{ A}_{n-2N,0}  \\
iv').&\ \bar{\c V}= p\c A^{(1)} = p\c D^{(1)}\hat{\c U} -p\c D^{(1)}\hat{\c T}_0 \underline{\eta}_{1} + \epsilon p^{n-2N} {\bold\Theta}^{(1)}(\underline{\eta}_{1, -\pi}^{\sigma \phi})^t(\hat{A}_{n-2N,0}\underline{\eta}_{1} - \hat{\c J}_{n-2N})  \\
v').&\ \bar{\c J}_{n-2N  } = p^{n-2N}\hat{\c J}_{n-2N} - p^{n-2N}\hat{A} _{n-2N,0}\underline{\eta}_{1}.         
\end{split}
\end{equation}

By (2.33) we see that the rows of $\c Y^{(1)}$ are (semi-linear) ``eigenvectors" of the semilinear map
\[
\c Z \rightarrow \c Z^{\sigma \phi} \c A^{(1)}
\]
with eigenvalues $\{1, p,\dots, p^{n-2N-2}\}$.
We will see below that $\c Y^{(1)}$ is a fundamental solution matrix of the differential equation (2.45) and therefore the rows of $\c Y^{(1)} $ form a basis for the solution space of (2.39) consisting of normalized period solutions. In particular, at $\hat{\lambda}_0 = 0$ the first row of $\c Y^{(1)}_{\hat{\lambda}_0}$, namely $\c Y^{(1)}_{\hat{\lambda}_0, 0}$, satisfies (dropping the subscript $\hat{\lambda}_0$)
\[ (\c Y^{(1)}_0)^{\sigma \phi} \c A^{(1)} = \c Y^{(1)}_0.\]
Then we will write explicitly $\c Y^{(1)}_0 = (\c Y^{(1)}_{0;(0,0)}, \c Y^{(1)}_{0;(0,1)}, \dots, 
 \c Y^{(1)}_{0, (0; n-2(N+1)} )$.
It is an immediate consequence of the fixed point relation (2.16) and the estimates for $\alpha_{i,j}$ with $j>i$ in (I) above that $|\underline{\eta}| < 1$. From the choice of initial data 
${\c Y}_{N;\hat{\lambda}_0}(0) = P_N(0) = I$ at $\hat{\lambda}_0 = 0$ it follows that $\underline{\eta} = (0,\dots,0)$.  It is useful to enumerate the rows and columns of $\c A^{(1)}$ as running from $0$ to 
$n-2N-2$ (instead of running from $1$ to $n-2N-1$). Let the entries of $\c A^{(1)}$ be denoted $\c A^{(1)} = (\alpha^{(1)}_{j,i})$ with $0 \leq i,j \leq n-2N-2.$  Let the entries of $\hat{\c A}^{(1)}$ be denoted $\hat{\c A}^{(1)} = (\hat{\alpha}^{(1)}_{j,i})$. It follows then from (2.40 iv') above that the following estimates hold for the entries of $\hat{\c A}^{(1)}$:
\begin{equation}
\begin{split}
i).\   &\c A^{(1)}(0) \  \text{is lower triangular with diagonal equal to} \  {\c D}^{(1)} = {\c D}_{N+1} = \{1, p, \dots,p^{n-2N-2}\}, \\
ii).\    &{\c A}^{(1)}(\lambda) \ \text{factors as follows:}  \  \c A^{(1)}(\lambda) = {\c D}^{(1)} \hat{\c A}^{(1)}(\lambda) \ \text{with} \ \hat{\c A}^{(1)}(\lambda)\text{ has integral coefficients, analytic on }|\lambda|\leq1, \\
iii).\  &\text{ord} \ \hat{\alpha}^{(1)}_{j,i} > 0 \ \text{for} \ i>j, \\
iv).\  &{ \hat\alpha}^{(1)}_{j,j}(\lambda) \ \text{is a principal unit for all} \ j \in \{0,1,\ldots,n-2N-2\}.
\end{split}
\end{equation}
The properties (ii--iv) hold for all  $|\lambda|\leq1$.         
Setting $\lambda = 0$ in (iv') above, and recalling that $\underline\eta(0)=(0,\dots,0)$, it follows that $\c A^{(1)}(0)=p^{-1}\c U(0)=\hat{\c U}(0)$. From (2.31), then (2.41)(i) above follows from the lower triangularity of ${\c A}(0)$ and $\bar{A}(0)$ given in (2.3) above. This discussion 
establishes (I) for $\c A^{(1)} (= \c A_{N+1}) $ in the proposed \textbf{K-Pack}$\{N+1\}$. Both conditions of II bis) for index $N+1$ have already been verified above. The first-row integrality and unit properties follow directly from Lemma~\ref{P}.

We turn next to the differential equation satisfied by the matrix of local functions $\c Y^{(1)}_{\hat{\lambda}_0}$. We note that if we set $\delta = \lambda\frac{\partial}{\partial \lambda}$ and
\begin{equation}
\bar{G} = \c Y_2 G \c {Y}_2^{-1} - (\delta \c Y_2)\c Y_2^{-1},
\end{equation}
then from the definition of $\c Y_1$ and (II) of the definition of a {\bf K-Pack\{$N$\}} above, it follows that $\c Y_1$ is a collection of local solutions to the differential system
\begin{equation}
\delta {\c Y} = {\c Y} \bar{G}.  
\end{equation}
Since $\c Y_1$ is lower triangular, it follows from (2.43) that $\bar{G}$ is as well and we may write
\begin{equation}
\bar{G} =
\begin{vmatrix}
\bar{G}_{0,0} & 0 & 0\\
\bar{G}_{1,0} & G^{(1)} & 0\\                               \bar{G}_{2,0} & \bar{G}_{2,1} & \bar{G}_{2,2}\\
\end{vmatrix} .                     
\end{equation}
All entries in $\bar G$ have integral coefficients and are analytic and bounded by $1$ on $|\lambda|\leq1$. It is immediate then that 
\begin{equation}
\delta {\c Y}^{(1)} = {\c Y}^{(1)} G^{(1)},    
\end{equation}
where $\bar{G}(0) = H_N$ and $G^{(1)}(0) = H_{N+1}$. Of course, the differential system (2.45) has a full set of logarithmic solutions, say, 
$\c Y_{N+1;0} = \lambda^{H_{N+1}} P_{N+1}(\lambda)$ with $P_{N+1} $ holomorphic in $D(0,1^{-})$. If we take $\c Y_{N+1, \lambda_0} = \c Y^{(1)}_{\lambda_0}$ at $\lambda_0=0$, then it follows easily  that $\c Y_{N+1, \lambda_0} = \lambda^{H_{N+1}} P_{N+1}$ is a fundamental system of solutions of (2.45) at $\lambda_0 = 0$ with initial condition $P_{N+1}(0) = I.$.
It follows at the origin from 
\[
{\c Y}^{\sigma \phi} {\c A} {\c Y}^{-1} = \c D_0
\]
for our given normalized solutions, that
\[
\c Y_1^{\sigma \phi} \bar{\c  A} = \c D_0 \c Y_{1}
\]
and from this
\begin{equation}
(\c Y^{(1)})^{\sigma \phi} \c A^{(1)} = \c A^{(1)}(0) \c Y^{(1)}.
\end{equation}

We complete the proof of the theorem by choosing as the 
{\bf K-Pack\{$N+1$\}} the collection of $n-2N -1$ by $n-2N-1$ matrices
\[
\{ \c A^{(1)}, G^{(1)}, {\Theta}_{N+1}, \c Y^{(1)} \}
\]
where ${\Theta}_{N+1}= -{\Theta}_N^{(1)}$ and for $\hat{\lambda}_0$ a unit in $\Gamma$ we take 
 $\c Y_{N, \hat{\lambda}_0}^{(1)}$ as the middle block of $\c Y_{1,\hat\lambda_0}$. On the other hand for $\hat{\lambda}_0 = 0$ we take the solution matrix $\c Y_{N+1,0} = \lambda^{H_{N+1}} P_{N+1}(\lambda)$ with a full set of logarithmic solutions.
 
We note from the definition in Corollary 2.3 above that $\c Y_2(0) = I_{n-2N +1}$, so that $\bar{A}_N(0) = M_0$ and in particular is lower triangular. It follows then that $\c A_{N+1}(0) (= \c A^{(1)}(0))$ is also lower triangular with diagonal entries $(1,p,\dots,p^{n-2N-2})$.
Also our choice of normalized solution matrix at $\hat{\lambda}_0 = 0$ means that its first row $\c Y^{0}_{N+1}(0) = (1,0,\dots,0)$ and the lead entry satisfies $\c F_{N+1}(0) = 1$. 

The symplectic relation (III) for $\c Y^{(1)}$, namely,
\[ 
\c Y^{(1)} {\Theta} \c Y^{(1)}_{-\pi} = \theta^{(1)},
\]
(with $\Theta = \Theta_{N+1}$ ) then follows by comparing the (1,1) blocks in the  symplectic relation for $\c Y_1$
\[
\c Y_1 {\Theta} \c Y_{1,-\pi}^t = \theta.
\]
This relation is itself a consequence of the definition of $\c Y_1$, and the symplectic relation (2.26) for $\c Y_2$.

In similar fashion, starting from the functional equation (IV) for the Frobenius matrix ${\c A}$, one obtains a functional equation for the matrix $\bar{A}$, namely
\[
\bar{A}{\Theta}_N \bar{A}_{-\pi}^t = p^{n-2N}{\Theta}_N,
\]
which is obtained by substituting $(\c Y_2^{-1})^{\sigma \phi} \bar{A} \c Y_2$ for $\c A$ (and the corresponding relation for $\c A_{-\pi}$) in the functional relation (IV). Comparing the middle blocks and dividing by $p^2$ gives the required functional equation for index $N+1$:
\[
\c A^{(1)}\Theta_{N+1}(\c A^{(1)}_{-\pi})^t=p^{n-2N-2}\Theta_{N+1}.
\]

\end{proof}

The following result then is immediate and gives $p$-adic analytic formulas for the eigenvalues in the case of the $L$-function associated with the $n$-dimensional hyperkloosterman sum (assuming $p > n+1$). 
\begin{theorem}
Let $n$ be a fixed positive integer.  Then for $0 \leq N \leq \left \lfloor{\frac{n}{2}}\right \rfloor$ let $\{ \c A_N, G_N, {\Theta}_N, \c Y_N \}$ be the {\bf K-Pack\{$N$\}} of normalized solution matrices of the preceding result,  let $\c F_N(\lambda)$ be the lead entry of the unique bounded solution $\c Y^{(0)}_N(\lambda)$ at the origin, normalized by $\c F_N(0)=1$, of the differential equation
\begin{equation}
\lambda\frac{\partial}{\partial \lambda} \c Y =  \c Y G_N(\lambda).
\end{equation}
Then for $\lambda_0 \in \mathbb{F}_q^{*}$, with Teichmuller lift $\hat{\lambda}_0 \in \Gamma$,   the collection $\{ \omega^{(N)}(\lambda_0)\}_{N=0}^{\left \lfloor{\frac{n}{2}}\right \rfloor}$, where $\omega^{(N)}(\lambda_0)= q^N  \frac{\c F_N(\lambda)}{(\c F_N)^{\sigma^a} (\lambda^q)}|_{\lambda = \hat{\lambda}_0}$,  is the collection of $\left \lfloor{\frac{n}{2}}\right \rfloor + 1$  eigenvalues of $\c A((\hat{\lambda}_0)^{q/p})\cdots\c A(\hat{\lambda}_0^p)\c A(\hat{\lambda}_0)$
having the smallest $ord_q$ values  running from $0$ to $\left \lfloor{\frac{n}{2}}\right \rfloor $. The eigenvalues of largest $ord_q$ values are then the collection 
$\{ \frac{q^{n }}{\omega^{N}_{-\pi}(\lambda_0)} \}_{N=0}^{\left \lfloor{\frac{n-1}{2}}\right \rfloor}$. In the case $n$ is even, we use the observation that $det(\c A_{n-2N}) = p^{\frac{n(n+1)}{2}}$ to conclude that the eigenvalue having $ord_q = \left \lfloor{\frac{n}{2}}\right \rfloor$ may also be  given by
\[
q^{\frac{n}{2}} \prod_{N=0}^{\frac{n}{2}-1} \frac{\omega^{(N)}_{-\pi}(\lambda_0)}{\omega^{(N)}(\lambda_0)}.
\]
\end{theorem}
We observe for future use the following generalization of Dwork (6.11) (Bessel Functions) which follows immediately from the above work. Let $n$ be fixed and $p > n+1$. We recall from Section $2$, we have that for the {\bf K-Pack\{$N$\}} the column vector $\underline{\zeta}_N  = (\zeta_{N;1}, \zeta_{N;2},\dots,\zeta_{N;n-2N-1})^t$ appears in the matrix $\c{Y}_{1}$ (Equation $(2.28)$).  Here 
$ \zeta_{N,i}U$ is the $(i,0)$ entry of the (local) matrix $\c Y_1$. More precisely in our earlier notation, let $\c{Y}_N^{(1)}(\lambda)$ denote the second row vector solution of the system $(2.1)$ in our normalized solution matrix $\c{Y}_N$ at the origin. In an obvious sense this is the log solution of this system and has lead entry $U_N(\lambda) \zeta_{N;1}(\lambda) = \c{F}_N(\lambda) \log(\lambda) + \c{G}_N(\lambda)$ (recall that at $\hat{\lambda}_0 = 0,  U_N = \c F_N$).For every $N$, $0  \leq N \leq \left \lfloor{\frac{n-1}{2}}\right \rfloor$, we have
\begin{equation}
p\zeta_{N;1} - \zeta_{N;1}^{\sigma \phi} \in p \c{O}_U[[\lambda]]. 
        \end{equation}
\section{Explicit formulas}

In the previous section, we produced $p$-adic analytic formulas for all the eigenvalues of Frobenius in the case of the $n$-variable Kloosterman sum. These formulas are expressed in terms of  solutions of the deformation equation and the differential equations obtained successively using the inductive process described in section 2. Here we give a description of these formulas in terms of classical functions that arise from the deformation equation and its exterior powers.  To that end we review linear algebra operations and the concomitant Frobenius action on the Dwork dual cohomology space and the consequent Frobenius action on these differential equations.

Let $\c L$ be a linear algebra operation, such as taking tensor  power, symmetric power, exterior power, or some composition of these operations. If $\c K$ denotes the Dwork dual space, we denote by $\c L \c K$ the space obtained from the Dwork dual space $\c K$ by acting on it by $\c L$. Our focus of interest here will be the case $\c L =\bigwedge^j$ , with $1\leq j\leq n+1$.  If we fix the  usual dual basis $\c B^{*} =\{ e^{*}_j \}_{j=0}^n $ for $\c K$,  it is convenient to choose the basis of $\bigwedge^j \c K$ as follows. Let $\c E_{n,j} = \{\underline{k} =(k_1, k_2, \cdots,  k_j) \in \mathbb Z^j| 0 \leq k_1 < k_2 < \cdots < k_j \leq n \}$. Let $ \bigwedge^j \c B^{*}  = \{\underline{e}^{*}_{\underline{k}} = e_{k_1}^{*} \wedge e_{k_2}^{*} \wedge \cdots \wedge e_{k_j}^{*} | \underline{k} \in \c E_{n,j}\}$ where we have   used $\underline{k}$ to denote the index. We define the weight of $\underline{k}$, denoted $|\underline{k}|$, to be $\sum_1^j k_i$. To fix ideas and to be as explicit as possible  we will choose a linear order for the basis $\bigwedge^j \c B^{*}$ as follows. Fix $\underline{e}^{*}_{\underline{k}} < \underline{e}^{*}_{\underline{\ell}}$ if $\text{weight} ({\underline{k}}) < \ \text{weight} (\underline{\ell})$.  If $\text{weight}(\underline{k}) = \text{weight}(\underline{\ell})$, we say $\underline{e}^{*}_{\underline{k}} < \underline{e}^{*}_{\underline{\ell}}$ provided $\underline{k} $ is less than $\underline{\ell}$ in the usual lexicographic ordering.
Let $\c A$ be the matrix of Frobenius acting on $\c K$ with respect to the basis $\c B^{*}$. We denote then the matrix of the $\bigwedge^j$-th power of Frobenius with respect to the ordered basis $\bigwedge^j \c B^{*}$ by $\bigwedge^j \c A=:\mathfrak{A}^{(j)}$  (omitting whenever possible explicit mention of the basis).  Then $\bigwedge^{(j)} \c K$ is again a module with connection, say $ \nabla_{\wedge^j}(\lambda\frac{\partial}{\partial \lambda})$. Of course, the action of this connection is governed as usual by the Leibniz rule acting on exterior products.  Let us fix ideas and notation so that if $G$ is the connection matrix of $\nabla$ with respect to the basis $ B^{*}$, then we will write $\bigwedge^j G = :\mathfrak{G}^{(j)}$ for the connection matrix of $\nabla_{\wedge^j}$ with respect to the ordered basis $\bigwedge^j \c B^{*}$. The horizontal sections $Y_{\wedge^j, \hat{\lambda}_0}$ in a neighborhood of the point $\hat{\lambda}_0$ for the connection $\nabla_{\wedge^j}$ satisfy the linear $\c N^{(j)}$ by $\c N^{(j)}$ system for the j-th exterior power with $\c N^{(j)}=\binom{n+1}{j}$.
\begin{equation}
\lambda\frac{\partial }{\partial \lambda} (Y_{\wedge^j, \hat{\lambda}_0}) =Y_{\wedge^j, \hat{\lambda}_0} \mathfrak{G}^{(j)}. 
\end{equation}
It is a very simple argument from the formalism of Dwork's deformation theory that on $\bigwedge^j \c K$, the $j$-th exterior power of the Dwork dual space, the Frobenius $\bigwedge^j \text{Frob}$ and the connection $\nabla_{\bigwedge^j}$ commute in the following sense. Let $\hat{\lambda}_0$ be a point of $\c O_{\c T}$ satisfying $\sigma(\hat{\lambda}_0) = \hat{\lambda}_0^p$. Let $ Y_{\wedge^j, \hat{\lambda}_0}$ be a solution matrix of 
(3.1), then
\begin{equation}
\ Y_{\wedge^j, \hat{\lambda}_0}^{\sigma}(\lambda^p) \wedge^j  \c A(\lambda) = \c M^{(j)} \ Y_{\wedge^j, \hat{\lambda}_0}(\lambda).  
\end{equation}
As usual, $\c M^{(j)}$ denotes a locally constant matrix.

The dual space $\c K$ is a free $\c C_0$-module of rank $n+1$ and the $j$-fold exterior power is free of rank $\binom{n+1}{j}$. A solution of equation (3.1) is a (row) vector of rank $\c N^{(j)}$, say $\underline{y}$, such that 
$\lambda \frac{\partial \underline{y}} {\partial \lambda}(\lambda)  = \underline{y}(\lambda) \mathfrak{G}^{(j)}.$
In particular, given $\underline{k}=(k_1,\ldots,k_j)\in \c E_{n,j}$ and $j$ solutions $\{ y_{k_1}, \dots , y_{k_j} \}$ of (2.1), then $\mathfrak{G}^{(j)}$ has been defined so that $\underline{y}_{\underline{k}} =y_{k_1} \wedge  y_{k_2} \wedge \cdots \wedge y_{k_j}$ is a solution of (3.1). 
At a nonzero Teichm\"uller center, we now consider the case where $y_i$ is the normalized solution of (2.1) associated with the (semi-linear) eigenvalue of Frobenius $p^i$. 
Then $\underline{y}_{\underline{k}}$ is again a normalized solution vector for (3.1) having (semilinear) eigenvalue $p^{|\underline{k}|}$ for $ \mathfrak{A}^{(j)}$ .

For notational convenience we denote the entries in the matrix $\mathfrak{A}$ as $\mathfrak{A}_{\underline{i},\underline{\ell}}$  with $\underline{i}$ and $\underline{\ell}$ in $\c E_{n,j}$ .

For example, in the case of  interest here, when $\c L= \bigwedge^j$ is the $j$-th exterior power, we proceed as follows.  
The entries of the matrix $\mathfrak{A}^{(j)}$ are the $j$ by $j$ minors of the matrix $\c A$ indexed by the $j$-th exterior products $\{e^\ast_{k_1} \wedge e^\ast_{k_2} \wedge\cdots\wedge e^\ast_{k_j} \mid 0 \leq k_1 < k_2 <\dots <k_j \leq n \}$ of $\bigwedge^j \c B^{*}$. The differential equation (3.1) is in this case the usual $j$-th exterior power of the differential equation (1.2) of the Introduction. 
One may obtain solutions to the differential system (3.1) starting  with a fundamental solution matrix, say $\c Y_{\hat{\lambda}_0}$ at $\hat{\lambda}_0$ for (1.2). A  fundamental solution matrix for (3.1), $\mathfrak{Y}^{(j)}$, may then  be obtained  by taking for its  $\underline{k}$-th row  $\mathfrak{Y}_{\underline{k}}=\sum_{\underline{j} \in \c E_{n,j}} \mathfrak{Y}_{\underline{k}, \underline{j}} e^{*}_{\underline{j}}$ and the coefficient $\mathfrak{Y}_{\underline{k}, \underline{j}}$ is the j by j minor with rows indexed by $\underline{k}$ and columns indexed by $\underline{j}$ of the chosen fundamental solution matrix $\c Y$. We examine these relations for $\hat{\lambda}_0 \in \Gamma$ in the particular cases when the chosen fundamental solution matrix $\c Y_{\hat{\lambda}_0}$ for (1.2) is the fundamental normalized period matrices of Section 2 for each $\lambda_0 \in \Gamma$ .

\begin{lemma}
Let $\c A$ denote the Frobenius matrix with respect to the basis $\c B^{*}$ as in Section  2. We assert the following estimates hold for the matrix $\mathfrak{A}^{(j)}$
of Frobenius acting on $\bigwedge^j \c K$ with respect to the basis $\bigwedge^j \c B^{*}$. We will use the notation $\underline{0}$ to denote the indices $\{0,1,2, \dots, j-1\}$.

 The matrix $\mathfrak{A} = \bigwedge^j \c A$ has entries $\mathfrak{A}_{\underline{k}, \underline{m}}$ where $\underline{k}$ and $\underline{m}$ lie in $\c E_{n,j}$. 
We claim the following estimates hold for the entries of the matrix $\frak A$.
\begin{align}
(i).\ &{\rm ord}_p (\mathfrak{A}_ {\underline{i},\underline{\ell}}) \geq |\underline{i}| \  \nonumber \\
(ii).\ &{\rm ord}_p( \mathfrak{A}_{\underline{i},\underline{i}}) = |\underline{i}|,  \\ 
(iii).\ &{\rm ord}_p (\mathfrak{A}_{\underline{i},\underline{\ell}}) > |\underline{i}|, \ \text{if} \ |\underline{i}| < |\underline{\ell}|, \ \text{or if} \ 
|\underline{i}| = |\underline{\ell}| \ \text{and} \ \underline{i} \neq \underline{\ell}. \nonumber
\end{align} 
Furthermore, $\mathfrak{A}_{\underline{i}, \underline{i}} = p^{|\underline{i}|} \hat{\mathfrak{A}}_{\underline{i}, \underline{i}}$ with 
$ \hat{\mathfrak{A}}_{\underline{i}, \underline{i}}$ a principal unit for all $\underline{i} \in \c E_{n,j}.$
\end{lemma}

\begin{proof}
Clearly if $\underline{i} = (i_0, i_1,\dots,i_{j-1})$ then the minor $\mathfrak{A}_{\underline{i}, \underline{\ell}}$ consists of sums of $j$-fold products from the matrix $\c A$ every one of which has an entry from each of the rows $i_{0}, \cdots i_{j-1}$. The first part then follows from the estimates (2.3)(ii) in section 2 above. Similarly, let us assume the hypothesis of (3.3)(iii), namely  either $|\underline{i}| < |\underline{\ell}|, $ or if $ 
|\underline{i}| = |\underline{\ell}| \ \text{and} \ \underline{i} \neq \underline{\ell}$. 
In this case, each summand in  $\mathfrak{A}_{\underline{i}, \underline{\ell}}$ up to a sign, has the form $\alpha_{i_0,\ell_0} \alpha_{i_1 ,\ell_1} \cdots \alpha_{i_{j-1}.\ell_{j-1}}$ with at least one factor $\alpha_{i_s, \ell_s}$ having second co-ordinate $\ell_s > i_s$.
This ensures by (2.3)(iii) that $\text{ord}_p \ (\alpha_{i_s,\ell_s}) > i_s$ and $ \text{ord}_p (\alpha_{i_0,\ell_0} \alpha_{i_1 ,\ell_1} \cdots \alpha_{i_{j-1}.\ell_{j-1}}) > \ |\underline{i}|$. 
This concludes the proof of (3.3)(iii). Precisely this same argument shows that every summand in $\mathfrak{A}_{\underline{i}, \underline{i}}$ other than the product $ \alpha_{i_0, i_0}\alpha_{i_1, i_1}\cdots\alpha _{i_{j-1}, i_{j-1}}$ also satisfies that its $\text{ord}_p$ is strictly greater than $\text{weight} \ (\underline{i})$ and (3.3)(ii) follows from (2.3)(iv).
\end{proof}

The matrix $\frak A$ is an $N^{(j)}$ by $ N^{(j)}$ matrix with entries in $\ C_0$ and $N^{(j)} = \binom{n+1}{j}$. When no confusion arises we will drop the superscript $j$ and write $N$ for $N^{(j)}$. We have described above a choice for the linear order on the basis $\bigwedge^j \c B^{*}$. Using this order our basis is $e^*_{\underline{i}_0}, \cdots,  e^*_{\underline{i}_{N-1}}$ .To simplify notation we use, abusively, the notation $\underline{s}$ to denote the basis element $e^*_{\underline{i}_{s}}$.  We have already noted that in the linear order on the basis, we write $\underline{0}$ for the first basis element $e^*_{\underline{i}_0} =e^{*}_0 \wedge e^{*}_1 \wedge \cdots e^{*}_{j-1}$.

The following corollary is an immediate consequence of Lemma 3.1.  The proof is immediate.
\begin{corollary}
The Frobenius matrix acting on $\bigwedge^j \c K$ with respect to the basis $\bigwedge^j \c B^{*}$ can be factored 
\[
\mathfrak{A} = \mathfrak{D} \hat{\mathfrak{A}}
\]
where $\mathfrak{D}$ is a diagonal $N$ by $N$ matrix with diagonal entries $\{ p^{|\underline{0}|}, p^{|\underline{1}|},\dots, p^{|\underline{N-1}|} \}$ where we have abusively denoted the ith element of $\c E_{n,j}$ by $\underline{i}$ in the chosen linear order on $\c E_{n,j}$. Furthermore, the following estimates hold:\\
(i). The entries $\hat{\mathfrak{A}}_{\underline{i}, \underline{\ell}}$ of $\hat{\mathfrak{A}}$ have integral coefficients and are analytic and bounded by $1$ on $D(0,1^+)$. \\
(ii). The diagonal entries of $\hat{\mathfrak{A}}$ are principal units. \\  
(iii). The entries $\hat{\mathfrak{A}}_{\underline{i}, \underline{\ell}}$ with $\underline{i} < \underline{\ell}$ are p-integral non-units.
\end{corollary}

We consider next the Frobenius structure on the $j$-th exterior power of our  confluent hypergeometric differential equation (1.2). For any $j$-fold collection of solutions of (1.2), say $\c Y_{\underline{k}} = \{\c Y^{k_1}, \c Y^{k_2},\dots,\c Y^{k_j} \}$ indexed by $\underline{k} = \{k_1, k_2,\dots,k_j\}$ with $\underline{k} \in \c E_{n,j}$ arranged according to the chosen linear order on $\wedge^j \c B$. We obtain thus a $1 \times N^{(j)}$  row vector solution  to the differential equation (3.1) with $\underline{k}$ coordinate equal to  $\mathfrak{Y}^{(j)}_{\underline{k}}=\c Y^{k_1} \wedge \c Y^{k_2} \wedge \cdots \c Y^{k_j}$. Recall that if $\hat{\lambda}_0 \in \Gamma $ then each solution $\c Y_{\ell}$ of (1.2) at $\lambda_0$ is a 1 by (n+1) row vector $ \c Y_{\ell} = (\c Y_{\ell,0}, \c Y_{\ell,1},\dots,\c Y_{\ell,n})$.  
(we are suppressing in the notation the subscript $\hat{\lambda_0}$). We form the solution vector $\mathfrak{Y}^{(j)}_{\underline{k}}$, an $N^{(j)}$-vector, of the $j$-th exterior power of the differential system (1.2), taking as coefficients all the $j$ by $j$ minors from the $j$ rows indexed by $\underline{k}$. Namely, let   
\[
\mathfrak{Y}_{\underline{k},\underline{m}} = \det(\mathfrak{Z}_{\underline{k},\underline{m}})
\]
where 
\[
\mathfrak{Z}_{\underline{k},\underline{m}} = 
\begin{vmatrix}
\c Y_{k_1,m_1} & \c Y_{k_1,m_2} &  \dots  & \c Y_{k_1,m_j} \\
\c Y_{k_2,m_1} & \c Y_{k_2,m_2} & \dots   & \c Y_{k_2,m_j} \\                   
\vdots & \vdots & \vdots & \vdots \\                                                   
\c Y_{k_j,m_1} & \c Y_{k_j,m_2} & \dots & \c Y_{k_j,m_j}  \\              
\end{vmatrix}.
\]

Then the element $\frak Y_{\underline{k}}(\lambda) = \sum_{\underline{m} \in \c E_{n,j}} \frak Y_{\underline{k}, \underline{m}}(\lambda) e^{*}_{\underline{m}}  =(\frak Y_{\underline{k},\underline{0}}(\lambda) , \cdots \frak Y_{\underline{k}, \underline{N-1}}(\lambda)), $  is a solution of the $j$-th exterior power differential equation 
\begin{equation}
\lambda\frac{\partial}{\partial \lambda} \mathfrak{Y} = \mathfrak{Y} \mathfrak{G}^{(j)}.   
\end{equation}
The Frobenius action on solutions then has the form (for a fundamental system of solutions $\mathfrak{Y}_{\hat{\lambda}_0}$) in a neighborhood of a point $\hat{\lambda}_0$ satisfying $\sigma(\hat{\lambda}_0)= \hat{\lambda}_0^p$ as follows:
\begin{equation}
\mathfrak{Y}_{\hat{\lambda}_0}^{\sigma \phi} \mathfrak{A}(\lambda) = \mathfrak{M}_{\hat{\lambda}_0} \mathfrak{Y}_{\hat{\lambda}_0}.
\end{equation}

As before, the matrix $\mathfrak{M}$ is locally constant, so constant  on the open unit disk centered at $\hat{\lambda}_0$. If we consider (3.4) when our starting point is the matrix $ \mathfrak{Y}$ of j-fold collections of normalized period solutions, then in the case $\hat{\lambda}_0 \in \Gamma, \hat{\lambda}_0 \neq 0$,  we have (3.5) with $\frak{M}_{\hat{\lambda}_0} = \frak{D}$,
the diagonal $N$ by $N$ matrix with entries $p^{|\underline{0}|}, \cdots p^{| \underline{N-1}|} $ . If we divide both sides of (3.5) by $p^{|\underline{0}|}$  
and denote $\bar{\mathfrak A}=p^{-{|\underline{0}|}} \mathfrak{A}$ and $\frak{\bar{D}} = p^{-{|\underline{0}|}} \frak{D}$,
then we have from Corollary 3.2, 
\begin{equation}
\frak{Y}^{\sigma \phi}\bar{\mathfrak{A}} = \mathfrak{\bar{D}} \frak{Y}.
\end{equation}
In particular the first row of (3.6) is at $\hat{\lambda}_{0} \in \Gamma\setminus\{0\}$
\begin{equation}
(\mathfrak{Y}^{\underline{0}})^{\sigma}(\lambda^p) \bar{\mathfrak{A}}(\lambda) =  \frak{Y}^{\underline{0}}(\lambda)  
\end{equation}
with eigenvalue $1$.
        An immediate consequence of (3.7) and Corollary 3.2 is that for every $\hat{\lambda}_0 \in \Gamma, \hat{\lambda}_0 \neq 0, \mathfrak{Y}^{\underline{0}}$ has $p$-integral coefficients in $D(\hat\lambda_0, 1^{-})$ of the series in $\lambda-\hat{\lambda}_0$ appearing in $\mathfrak{Y}^{\underline{0}}$ and has first entry $\mathfrak{Y}^{{\underline{0}}}_{\hat{\lambda}_0;\underline{0}}(\lambda)$ a $p$-adic unit for all $\lambda \in D(\hat{\lambda}_0, 1^{-})$ as in proof of Lemma 2.2.
        
For $\hat{\lambda}_0 = 0$, we obtain a similar result making use of the logarithmic structure for the normalized solution of (1.2) at the origin. We are particularly interested in the behavior of solutions in the neighborhood of the singular point $\hat{\lambda}_0 = 0$ since the differential system (3.4) and its solutions there make classical sense (working say over $\mathbb{C}$) and the solutions there have both $p$-adic meaning (carrying arithmetic import) and a classical  interpretation over $\mathbb{C}$ as well. 
We have then at the origin the following.

\begin{lemma}
We will write $\frak{F}^{(j)}_0(\lambda)$ for the function $\frak{Y}^{(j)}_{0;\underline{0},\underline{0}}(\lambda)$ which is the determinant of the j by j minor of the solution matrix  at $\hat{\lambda}_0 = 0$ with rows and columns indexed by $\underline{0} = (0, 1,\dots, j-1)$. This is the ``first" entry of the ``first" 1 by $N^{(j)}$-vector solution of $(3.4)$ namely $\frak{Y}^{(j)}_{0;\underline{0},\underline{0}}(\lambda)$ with row and column indexed by $\underline{0}$ at the origin. More precisely, we have 
\[
\mathfrak{Y}^{(j)}_{0;\underline{0},\underline{0}}(\lambda) = \det(\lambda^{H^{(j)}} P^{(j)}(\lambda))
\]
where $H^{(j)}$  (respectively $P^{(j)}$) is the $j \ \text{by} \ j$  submatrix  of $H$ (respectively $P$) whose rows and columns are indexed by $\underline{0}=(0,1,2,\dots,j-1)$. This is the first entry of this 1 X $N^{(j)}$-vector solution of the differential system  $(4.4)_j$. Since $\det(\lambda^{H^{(j)}}) = 1$, $\mathfrak{Y}^{(j)}_{0; \{\underline{0}, \underline{0} \} }$ is analytic in $D(0,1^{-})$, as are the other entries of this ``first" row (with row index $\underline{0}$) of the solution vector $\mathfrak{Y}^{(j)}_{0;\underline{0}}$. Furthermore, the series entries (in $\lambda$) of this row have p-integral coefficients so are convergent and bounded by $1$ in $D(0,1^{-})$ and $\frak{Y}^{(j)}_{0:\{\underline{0}, \underline{0} \}}$ is a unit there. 
\end{lemma}

\begin{proof}
The first $j$ rows of $\lambda^H P$ involve only the first $j$ rows of $P$, with left factor $\lambda^{H^{(j)}}$. Its determinant is $1$, so every coordinate of $\mathfrak Y^{(j)}_{0;\underline0}$ is a minor of these rows of $P$ and converges in $D(0,1^-)$. Since $P(0)=I$, this row has initial value $(1,0,\ldots,0)$. At the origin the constant Frobenius matrix is $\bigwedge^j\c D_0$. Because $\c D_0$ is lower triangular, its first exterior row has only its first entry nonzero, equal to $p^{|\underline0|}$. Consequently $(\mathfrak Y^{(j)}_{0;\underline0})^{\sigma\phi}\bar{\mathfrak A}=\mathfrak Y^{(j)}_{0;\underline0}$. The coefficient argument in Lemma \ref{P}, applied to this row and $\bar{\mathfrak A}\equiv E_{00}\pmod\pi$, proves that its first coordinate lies in $1+\pi\lambda\c O_{\c U}[[\lambda]]$ and its other coordinates in $\pi\lambda\c O_{\c U}[[\lambda]]$. This proves the assertions. Moreover, the uniqueness of $P$ with $P(0)=I$ for the differential equation with coefficients in $\Omega_1$ shows that these origin series have coefficients in $\Omega_1$ and are fixed by $\sigma$.

\end{proof}

Now block-decompose the Frobenius matrix $\mathfrak{A}$ into $\mathfrak{A} = p^{|\underline{0}|} \bar{\mathfrak{A}}$ and block divide as follows:
\[
\bar{\mathfrak{A}} = 
\begin{vmatrix}
\bar{\mathfrak{A}}_1 & \bar{\mathfrak{A}}_2  \\
\bar{\mathfrak{A}}_3  & \bar{\mathfrak{A}}_4  \\
\end{vmatrix},
\]
where $\bar{\mathfrak{A}}_1$ is the 1 by 1 matrix with entry $p^{-|\underline0|}\mathfrak{A}_{\underline{0}, \underline{0}}$, a $p$-adic unit, and $\bar{\mathfrak{A}}_2$ is a 1 by $N^{(j)}-1$ row vector with $p$-integral non-unit entries for $\lambda \in D(0,1^{-})$ and $\bar{\mathfrak{A}}_3$ (respectively $\bar{\mathfrak{A}}_4$) is an $N^{(j)}-1$ by 1 column matrix (respectively an $N^{(j)}-1$ by $N^{(j)}-1$) matrix with entries which are $p$-integral non-units for $\lambda \in D(0,1^{-})$). It follows then from (3.7) that $\mathfrak{Y}_{\hat{\lambda}_0,\underline{0}}$ for each $\hat{\lambda}_0 \in \Gamma$ (including $\hat{\lambda}_0 = 0$) is a fixed point of this Frobenius action when our original solution  matrix is the normalized solution matrix $\c Y_{\hat{\lambda}_0}(\lambda)$ as in the previous section. The contraction mapping argument of section 2 then establishes the following result:

\begin{theorem}
At the origin $\hat{\lambda}_0=0$, the function
\[
\frac{\mathfrak{Y}^{(j)}_{0;\underline{0},\underline{0}}(\lambda)}
{(\mathfrak{Y}^{(j)}_{0;\underline{0},\underline{0}})^{\sigma}(\lambda^p)}
\]
has analytic continuation to $D(0,1^+)$ and is a unit there.
If $\lambda_0\in\mathbb{F}_q^*$, with $q=p^a$, and
$\hat{\lambda}_0$ is the Teichm\"uller lift of $\lambda_0$, then
the unique unit eigenvalue $\bar{\mathfrak{U}}^{(j)}(\lambda_0)$
of the modified Frobenius
$\bar{\mathfrak{A}}(\hat{\lambda}_0^{p^{a-1}})
\bar{\mathfrak{A}}(\hat{\lambda}_0^{p^{a-2}})
\cdots\bar{\mathfrak{A}}(\hat{\lambda}_0)$
is given by the formula
\[
\bar{\mathfrak{U}}^{(j)}(\lambda_0)
=
\frac{\mathfrak{Y}^{(j)}_{0;\underline{0},\underline{0}}(\lambda)}
{(\mathfrak{Y}^{(j)}_{0;\underline{0},\underline{0}})^{\sigma^a}(\lambda^q)}
\bigg|_{\lambda=\hat{\lambda}_0}.
\]
The quotient is evaluated at $\hat{\lambda}_0$ after analytic
continuation. Furthermore, the eigenvalue
$\mathfrak{U}^{(j)}(\lambda_0)$ of the $q$-Frobenius matrix
$\mathfrak{A}(\hat{\lambda}_0^{p^{a-1}})
\cdots\mathfrak{A}(\hat{\lambda}_0)$
is then given by
\[
\mathfrak{U}^{(j)}(\lambda_0)
=
q^{|\underline{0}|}\bar{\mathfrak{U}}^{(j)}(\lambda_0).
\]
\end{theorem}

\begin{proof}
It is a consequence of the contraction mapping argument that 
\[
\bar{\mathfrak{U}}^{(j)}(\lambda) = \frac{\mathfrak{Y}^{(j)}_{0;  \{ \underline{0}, \underline{0} \}}(\lambda)}{(\mathfrak{Y}^{(j)}_{0; \{ \underline{0}, \underline{0} \}})^{\sigma^a}(\lambda^q)}
\]
is an analytic function on $D(0,1^{+})$ and for any $\hat{\lambda}_0 \in \Gamma$ such that $\sigma^a(\hat\lambda_0)=\hat\lambda_0$, $\bar{\mathfrak{U}}^{(j)}(\hat{\lambda}_0)$ is the unique unit eigenvalue of the Frobenius matrix $\bar{\mathfrak{A}}(\hat{\lambda}^{p^{a-1}}_0)\bar{\mathfrak{A}}(\hat{\lambda}^{p^{a-2}}_{0})\cdots\bar{\mathfrak{A}}(\hat{\lambda}_0)$ where the  matrix $\bar{\mathfrak{A}}^{(j)}(\lambda)$ acts on local solutions of $(3.4)_j$ (at $\hat{\lambda}_0$)
via
\[
\mathfrak{Y}^{(j)}(\lambda) \rightarrow (\mathfrak{Y}^{(j)})^{\sigma, \phi}(\lambda) \bar{\mathfrak{A}}^{(j)}(\lambda).
\]
 It follows then that $\mathfrak{U}^{(j)}(\hat{\lambda}_0) = q^{|\underline{0}|}\bar{\mathfrak{U}}^{(j)}(\hat{\lambda}_0)$ is the unique eigenvalue of $\mathfrak{A}(\hat{\lambda}^{p^{a-1}}_0)\mathfrak{A}(\hat{\lambda}^{p^{a-2}}_{0})\cdots\mathfrak{A}(\hat{\lambda}_0)$ having smallest $\text{ord}_q$ equal to $|\underline{0}|$. 
\end{proof} 

Finally, fix $\lambda_0 \in \mathbb{F}_q^{*}$.  We consider the $n$-dimensional multiple Kloosterman sums $Kl_{n+1}(x, \lambda_0)$ and the associated $L$-function $L(Kl_{n+1}, \mathbb{G}_{m,\mathbb{F}_q}, \lambda_0,T)^{(-1)^{n+1}} = \prod_{\ell=0}^n(1-T\omega_{\lambda_0.\ell})$, where we have denoted by $\omega_{\lambda_0,\ell}$ the unique reciprocal zero of the $L$-function  associated to the hyperkloosterman sum having $\text{ord}_q = \ell$. Of course, the set of reciprocal zeros of $L$ is the set of  the eigenvalues of the Frobenius matrix $\c A^{(q)}(\hat{\lambda}_0) = \c A(\hat\lambda_0^{p^{a-1}})\c A(\hat\lambda_0^{p^{a-2}})\cdots\c A(\hat{\lambda}_0)$ where $\hat{\lambda}_0$ is the Teichm\"uller lifting of $\lambda_0$ and $\c A$ is the Frobenius matrix as in $\textbf{K-Pack\{$0$\}}$ in section 2. Furthermore we know from the references that each $\omega_{\lambda_0.\ell}$ may be written as $q^{\ell} \bar{\omega}_{\lambda_0.\ell} $ with $\bar{\omega}_{\lambda_0.\ell} $ a principal $p$-adic unit. It has been the goal in this section to find a $p$-adic analytic function (given in terms of classical analytic functions) convergent on $D(0,1^{+})$ which interpolates the values $\{ \bar{\omega}_{\lambda_0.\ell} \}$ at the corresponding Teichm\"uller units $\{\hat{\lambda}_0 \}$. Let $\{ \eta^{*}_{\hat{\lambda}_0,\ell} \}$ be a basis of eigenvectors for $\c A^{(q)}$
\[
\eta^{*}_{\hat{\lambda}_0,\ell} \c A^{(q)}(\hat{\lambda}_0) = q^{\ell} \bar{\omega}_{{\lambda}_0.\ell} \eta^{*}_{\hat{\lambda}_0,\ell}.
\]
Let $\underline{k} = (k_1,\dots,k_j)$ and  denote $\eta^{*}_{\underline{k}} (= \eta^{*}_{k_1} \bigwedge\eta^{*}_{k_2}\bigwedge\cdots\bigwedge \eta^{*}_{k_j})$. It is the  $N^{(j)}$  vector with entries $(\eta^{*}_{\underline{k},\underline{0}}, \eta^{*}_{\underline{k},\underline{1}},\dots,\eta^{*}_{\underline{k},\underline{N^{(j)}-1}})$ where
\[
\eta^{*}_{\underline{k},\underline{m}}=
\begin{vmatrix}
\eta^{*}_{k_1,m_1} & \eta^{*}_{k_1,m_2} &  \dots  & \eta^{*}_{k_1,m_j} \\
\eta^{*}_{k_2,m_1} & \eta^{*}_{k_2,m_2} & \dots   & \eta^{*}_{k_2,m_j} \\
\vdots & \vdots & \vdots & \vdots \\                                                  
\eta^{*}_{k_j,m_1} & \eta^{*}_{k_j,m_2} & \dots & \eta^{*}_{k_j,m_j}  \\              
\end{vmatrix},
\]
and where the point $\hat{\lambda}_0$ is understood and removed from the notation for convenience.
Then $\eta^{*}_{\hat{\lambda}_0,\underline{k}} $ is an eigenvector of $\bigwedge^j \c A^{(q)}(\hat{\lambda}_0)$ associated with the eigenvalue $q^{|\underline{k}|} \prod_{\ell = 1}^{j} \bar{\omega}_{{\lambda}_0,k_{\ell}}$. But then the unique unit eigenvalue of  $q^{-|\underline{0}|}\bigwedge^j {\c A}^{(q)}(\hat{\lambda}_0)$ is
\[
\prod_{\ell =0}^{j-1} \bar{\omega}_{{\lambda}_0.\ell}.
\]
As we have shown above this proves $\bar{\frak{U}}({\lambda}_0) = \prod_{\ell =0}^{j-1} \bar{\omega}_{{\lambda}_0.\ell}$ and gives the following main result.
\begin{theorem}
With $\bar{\frak U}^{(0)}(\lambda)=1$, for $1\leq j\leq n+1$ we have
\[
\frak{U}^{(j)}(\hat{\lambda}_0) = \prod_{\ell = 0}^{j-1} \omega_{\bar{\lambda}_0, \ell} = q^{|\underline{0}|} \prod_{\ell = 0}^{j-1} \bar{\omega}_{\bar{\lambda}_0,\ell},
\]
so that, for $0\le j\le n$
\[
\omega_{\hat{\lambda}_0,j}  = q^j\frac{\bar{\frak{U}}^{(j+1)}(\lambda)}{\bar{\frak{U}}^{(j)}(\lambda)}\bigg|_{\lambda = \hat{\lambda}_0}. 
\]
For $j=n+1$, $\frak F_0^{(n+1)}=\det P=1$, because $\text{tr }G=0$, and hence 
$\bar{\frak{U}}^{(n+1)}=1.$ 
This expression then gives a $p$-adic analytic formula for the eigenvalues of Frobenius in terms of classical analytic functions coming from the exterior powers of the confluent hypergeometric differential equation underlying the Kloosterman sum. 
\end{theorem}

It may be useful to illustrate the case $j=2$ when $n \geq 3$. 
We consider the normalized solutions at $\lambda_{0} = 0$ associated with the eigenvalues having $\text{ord}_p = 0$ and $\text{ord}_p =1$ and higher. For simplicity we write 
$\c Y^{(0)} = (\c F_0, \frac{\lambda \partial}{\partial \lambda} \c F_0, \cdots (\frac{\lambda \partial}{\partial \lambda})^{n} \c F_0)$ and $\c Y^{(1)} = (\c F_1, \frac{\lambda \partial}{\partial \lambda} \c F_1, \cdots (\frac{\lambda \partial}{\partial \lambda})^{n} \c F_1) $ where $\c F_1 = \c F_0 \text{log} \ \lambda + \c G_0$. Then if $\frak{Y}^{(\underline{0})} = \c Y^{(0)} \wedge \c Y^{(1)}$, we have $\frak{U}^{(2)}(\hat{\lambda}_0) = \omega_{\bar{\lambda}_0,0} \omega_{\bar{\lambda}_0, 1} = q \prod_{\ell = 0}^{1} \bar{\omega}_{\bar{\lambda}_0,\ell},$
so that

i).    $\omega_{{\lambda}_0,0} = \bar{\mathfrak{U}}^{(1)}(\lambda) \bigg |_{\lambda =\hat{\lambda}_0} = \frac{\mathfrak{\c F}_0(\lambda)}{\mathfrak{\c F}_0(\lambda^q)} \bigg |_{\lambda =\hat{\lambda}_0}$,

ii). $\omega_{{\lambda_0}, 1}  = 
q\frac{\bar{\frak{U}}^{(2)}(\lambda)}{\bar{\frak{U}}^{(1)}(\lambda)}\bigg|_{\lambda = \hat{\lambda}_0}= q \frac{\c H (\lambda)}{\c H (\lambda^q)} \bigg |_{\lambda =\hat{\lambda}_0}$,

where $\c H= \c F_0 + \lambda\frac{\partial}{\partial \lambda}(\c G_0)-\c G_0 (\c F_0)^{-1} \lambda\frac{\partial}{\partial \lambda}(\c F_0) $.
One obtains similar,although more complicated, formulas for the other eigenvalues of larger $\text{ord}_p $.

\begin{center}
{\sc Acknowledgements}
\end{center}
H. Zhang is supported by National Key R\&D Program of China 2024YFA1014600, Jiangsu Provincial Scientific Research Center of Applied Mathematics (Grant No. BK20233002), National Natural Science Foundation of China (Grant No. 12301012), Natural Science Foundation of Jiangsu Province (Grant No. BK20230802),  and also by Zhishan Scholars Programs of Southeast University (Grant No. 2242025RCB0011).

\end{document}